\documentclass{elsarticle}

\usepackage{hyperref}

\usepackage{graphicx}
\usepackage{algpseudocode}

\usepackage{amsmath}
\usepackage{amsfonts}
\usepackage{amssymb}
\usepackage{bbm}

\usepackage{epstopdf}

\journal{{\phantom{none}}}

\newcommand{\ii}{\mathrm{i}}
\newcommand{\ee}{\mathrm{e}}

\newcommand{\lamo}{\overline\lambda}
\newcommand{\lamu}{\underline\lambda}

\newcommand{\ghost}{\textsf{GHOST}}

\begin{document}

\begin{frontmatter}

\title{High-performance implementation of Chebyshev filter diagonalization for interior eigenvalue computations}

 \author[address_greifswald]{Andreas Pieper}
 \author[address_erlangen]{Moritz Kreutzer}

 \author[address_greifswald]{Andreas Alvermann\corref{cor}}
 \ead{alvermann@physik.uni-greifswald.de}
 \cortext[cor]{Corresponding author}
 
 \author[address_wupp]{Martin Galgon}

 \author[address_greifswald]{Holger Fehske}
 \author[address_erlangen]{Georg Hager}
 \author[address_wupp]{Bruno Lang}
 \author[address_erlangen]{Gerhard Wellein}

 \address[address_greifswald]{Ernst-Moritz-Arndt-Universit\"at Greifswald, Germany}
 \address[address_erlangen]{Friedrich-Alexander-Universit\"at Erlangen-N\"urnberg, Germany}
 \address[address_wupp]{Bergische Universit\"at Wuppertal, Germany}

\begin{abstract}
We study Chebyshev filter diagonalization as a tool for the computation 
of many interior eigenvalues of very large sparse symmetric matrices.
In this technique the subspace projection onto the target space of wanted eigenvectors is approximated with filter polynomials obtained from Chebyshev expansions of window functions.
After the discussion of the conceptual foundations of Chebyshev filter diagonalization we 
analyze the impact of the choice of the damping kernel, search space size, and filter polynomial degree on the computational accuracy and effort, before we describe the necessary steps towards a parallel high-performance implementation.
Because Chebyshev filter diagonalization avoids the need for matrix inversion it can deal with matrices and problem sizes that are presently not accessible with rational function methods based on direct or iterative linear solvers.
To demonstrate the potential of Chebyshev filter diagonalization for large-scale problems of this kind we include as an example the computation of the $10^2$ innermost eigenpairs of a topological insulator matrix with dimension $10^9$ derived from quantum physics applications. 
\end{abstract}


\end{frontmatter}

\section{Introduction}

Computation of many interior eigenvalues and eigenvectors of a large sparse symmetric (or Hermitian) matrix is a frequent task required for electronic structure calculations in quantum chemistry, material science, and physics.
Starting from the time-independent Schr\"odinger equation $H \psi_i = E_i \psi_i$,
where $H$ is the Hamilton operator of 
electrons in the material under consideration,
one is interested in all solutions of this equation with energies $E_i$ within a small window inside of the spectrum of $H$.
Interior eigenvalues are required because electrical or optical material properties are largely determined by the electronic states in the vicinity of the Fermi energy, which separates the occupied from the unoccupied electron states.
The width of the energy window is given by, e.g., the temperature or the magnitude of an applied electric voltage.

Recently, the computation of interior eigenvalues has gained new thrust from the fabrication of graphene~\cite{CGPNG09} and topological insulator~\cite{Hasan10} materials, which exhibit unconventional electron states that differ from those observed in traditional semiconductors or metals.
In graphene, electrons at the Fermi level follow a linear dispersion reminiscent of relativistic Dirac fermions, while the  topologically protected surface states of topological insulators are natural candidates for quantum computing and quantum information storage.
Because the properties of graphene and topological insulator devices depend directly on such electron states close to the center of the spectrum, their modeling and design requires the solution of ever larger interior eigenvalue problems~\cite{SF12,SFFV12}.

Motivated by the rapidly developing demands on the application side we introduce in this paper a high-performance implementation of Chebyshev filter diagonalization (ChebFD), a straightforward scheme for interior eigenvalue computations based on polynomial filter functions.
ChebFD has a lot in common with the Kernel Polynomial Method (KPM)~\cite{WWAF06} and is equally well suited for large-scale computations~\cite{KHWPAF15}.
Our goal is to show that a well-implemented ChebFD algorithm can be turned into a powerful tool for the computation of many interior eigenvalues of very large sparse matrices.
Here we include results for the computation of $\simeq 10^2$ central eigenvalues of a $10^9$-dimensional topological insulator matrix on $512$ nodes of the federal supercomputer system ``SuperMUC'' at LRZ Garching, at $40$ Tflop/s sustained performance.
ChebFD is, to our knowledge, presently the only approach that can solve interior eigenvalue problems at this scale.

While the present paper focuses on high-performance applications of ChebFD,
and less on algorithmic improvements, we will discuss its theoretical foundations in considerable detail. 
We will analyze especially how the algorithmic efficiency of polynomial filtering depends on the number of search vectors in the computation, and not only on the quality of the filter polynomials.
Using a larger search space and lower filter polynomial degree can improve convergence significantly.
Conceptually, ChebFD is a massively parallel block algorithm with corresponding demands on computational resources, which makes performance engineering mandatory.

One should note that the use of Chebyshev filter polynomials has a long tradition in the solution of linear systems~\cite{Man77,Man78}, as a restarting technique for (Krylov) subspace iterations~\cite{Saad84,So02,Saad11}, and for the computation of extreme and interior eigenvalues in quantum chemistry applications~\cite{ZSTC06,Ne90,MT97a,MT97b}.
The latter methods can be divided roughly into two classes:
Methods which are based---explicitly or implicitly---on time propagation of a wave function with Chebyshev techniques and (re-)construction of eigenstates from Fourier inversion of the computed time signal (as, e.g., in Ref.~\cite{Ne90}),
and methods that use Chebyshev filters for iterated subspace projection steps (see, e.g., Ref.~\cite{dNPS14} for a recent contribution). ChebFD belongs to the latter class of methods.

The rationale behind the use of polynomial filters in ChebFD instead of rational filter functions (as, e.g., in FEAST~\cite{Pol09,2013-KraemerDiNapoliEtAl-DissectingTheFEAST-JComputApplMath:244:1-9} or the CIRR method~\cite{SS03,ST07}) or of spectral transformations of the eigenvalue problem~\cite{SBR08} is that sparse matrix factorization is not feasible
for the very large matrices in our target applications.
This prevents the use of direct solvers such as PARDISO~\cite{PARDISO} or ILUPACK~\cite{ILUPACK}, which have been used successfully for, e.g., Anderson localization~\cite{SBR08}.
In some applications the matrix is not even stored explicitly but constructed `on-the-fly'~\cite{Lin90,ALF11}. Therefore, we assume that the only feasible sparse matrix operation in ChebFD is sparse matrix vector multiplication (spMVM).

If a matrix $H$ is accessed exclusively through spMVMs
the only expressions that can be evaluated effectively are polynomials $p[H] \vec v$ of $H$ applied to a vector $\vec v$. Note that in such a situation also the use of iterative linear solvers would result in a polynomial expression. Therefore, we design and analyze the ChebFD algorithm entirely in terms of polynomial filters (see also Refs.~\cite{2014-Kraemer-IntegrationBasedSolvers-PhD,2015-GalgonKraemerLang-AdaptiveChoiceOf-PREPRINT, GKLAFPHKSWBRT16} for recent work on polynomial filters).

The paper is organized as follows.
In Sec.~\ref{sec:Poly} we describe the construction of Chebyshev filter polynomials, analyze the convergence of polynomial filtering as used in ChebFD, and obtain practical rules for the choice of the free algorithmic parameters.
In Sec.~\ref{sec:ChebFD} we introduce the basic ChebFD scheme and test the algorithm with a few numerical experiments.
In Sec.~\ref{sec:Impl} we describe the main performance engineering steps for our \ghost~\cite{GHOST}-based ChebFD implementation, which we use for the large-scale application studies in Sec.~\ref{sec:Large}.
We conclude in Sec.~\ref{sec:Conclude}, where we also discuss possible algorithmic refinements of the basic ChebFD scheme.

\section{Polynomial filter functions}
\label{sec:Poly}

\begin{figure}
\centering
\includegraphics[width=0.5\textwidth]{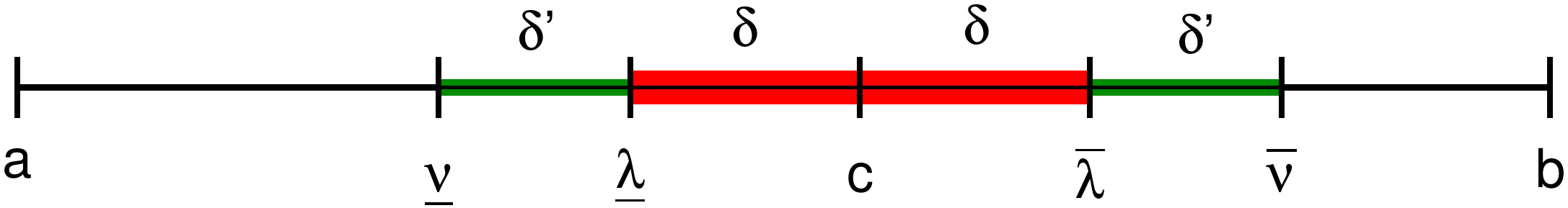}
\caption{Sketch of the interval configuration for interior eigenvalue computations:
The target interval $I_T = [\lamu,\lamo]$, with half width $\delta = (\lamo-\lamu)/2$, is contained within the search interval $I_S = [\underline\nu,\overline\nu]$, which is larger by a margin $\delta'$.
Both lie within an interval $[a,b]$ that encloses the spectrum of the matrix under consideration.
Typically, it is $\delta \ll b-a$ and the target interval center $c$ lies far away from $a$, $b$.
The symmetric situation $c \approx (a+b)/2$ is especially challenging.
}
\label{fig:sketch}
\end{figure}

The eigenvalue problem described in the introduction can be specified as follows:
Given a symmetric (or Hermitian) matrix $H$, compute the $N_T$ eigenpairs $(\lambda_i, \vec v_i)$ within a prescribed target interval $\lambda_i \in I_T = [\lamu,\lamo]$ (see Fig.~\ref{fig:sketch}).
Numerical convergence is achieved if the residual of every eigenpair satisfies $\| H \vec v_i - \lambda_i \vec v_i \| \le \epsilon$ for a specified accuracy goal $\epsilon$, and all eigenpairs from the target interval are found.
In current applications the target interval typically contains up to a few hundred eigenvalues of $H$,
which itself has dimension $10^6, \dots, 10^{10}$.
The precise value of $N_T$ is normally not known prior to the computation,
but can be estimated from the density of states (see Sec.~\ref{sec:Optimal}).

The optimal filter function for the above problem is the rectangular window function $W(x) = \Theta(x-\lamu) \Theta(\lamo -x)$, with the  Heaviside function $\Theta(x)$.
With this choice $W[H]$ is the exact projection onto the target space, that is,
$W[H] \vec v_i = \vec v_i$ for $\lambda_i \in I_T$ and
$W[H] \vec v_i = 0$ for $\lambda_i \not\in I_T$.
The eigenpairs in the target interval can then be obtained as follows:
Start with sufficiently many (random) search vectors $\vec x_k$, for $k=1, \dots, N_S$ with $N_S \ge N_T$,
and compute the filtered vectors $ \vec y_k =  W[H] \vec x_k $.
The vectors $\vec y_k$ lie in the target space,
but if $N_S > N_T$ they are linearly dependent.
Orthogonalization of the $\vec y_k$ with a rank-revealing technique gives a basis\footnote{It can happen that basis vectors are missing, but the probability of this event is negligible for $N_S \gg N_T$. In the ChebFD scheme of Sec.~\ref{sec:ChebFD} the probability is reduced even further by adding new vectors after each orthogonalization step.} of the target space.
The wanted eigenpairs are the Rayleigh-Ritz pairs of $H$ in this subspace basis.

\subsection{Construction of polynomial filter functions}

 Of course, a rectangular window function $W(x)$ cannot be represented exactly by a polynomial $p(x)$.
 The first step in polynomial filtering thus is the construction of a good polynomial approximation to $W(x)$. Theoretical arguments~\cite{Saad11} and numerical evidence (see Sec.~\ref{sec:ChebFD}) suggest that convergence improves with high-order polynomial approximations.
Since general optimization or fitting techniques become difficult to use for polynomials of high degree,
we construct the filter polynomials from Chebyshev expansions and kernel polynomials. This is motivated by our previous experience with the KPM~\cite{WWAF06}.

The construction of filter polynomials is based on the Chebyshev expansion
\begin{equation}\label{W1}
 W(x) = \sum_{n=0}^\infty c_n T_n(\alpha x + \beta) \;,  \qquad \text{ with } \alpha = \frac{2}{b-a} \;,  \;\; \beta = \frac{a+b}{a-b} \;,
\end{equation}
of $W(x)$ on an interval $[a,b]$ that contains the spectrum of $H$ (see Fig.~\ref{fig:sketch}).
Here, $T_n(x)= \cos (n \arccos x)$ is the $n$-th Chebyshev polynomial of the first kind, with recurrence
\begin{equation}
 T_0(x) = 1 \;, \quad T_1(x) = x \;, \quad T_{n+1}(x) = 2 x T_n(x) - T_{n-1}(x) \;.
\end{equation}
Notice the scaling of the argument $\alpha x + \beta \in [-1,1]$ of $T_n(\cdot)$ for $x \in [a,b]$.
The expansion coefficients $c_n$ are given by
\begin{align}
 c_n &= \frac{2 - \delta_{n0}}{\pi} \int_{-1}^1 W\Big( \frac{b-a}{2} x + \frac{a+b}{2} \Big)  T_n(x) (1-x^2)^{-1/2} \, dx \\
  & = \begin{cases} \frac1\pi \big( \arccos (\alpha \lamo + \beta ) -   \arccos (\alpha \lamu + \beta )  \big) \qquad & \text{ if } n = 0 \;, \\[0.5ex]
  \frac2 {\pi n} \Big( \sin \big(n \arccos (\alpha \lamo + \beta) \big) - \sin \big( n \arccos (\alpha \lamu + \beta ) \big)  \Big)  & \text{ if } n \ge 1 \;.  
  \end{cases}
\end{align}
Here, the first expression results from orthogonality of the $T_n(x)$ with respect to the weight function $(1-x^2)^{-1/2}$, the second expression is obtained by evaluating the integral for the rectangular function $W(x) = \Theta(x-\lamu) \Theta(\lamo -x)$.

\begin{figure}[h]
\hspace*{\fill}
\includegraphics[width=0.4\textwidth]{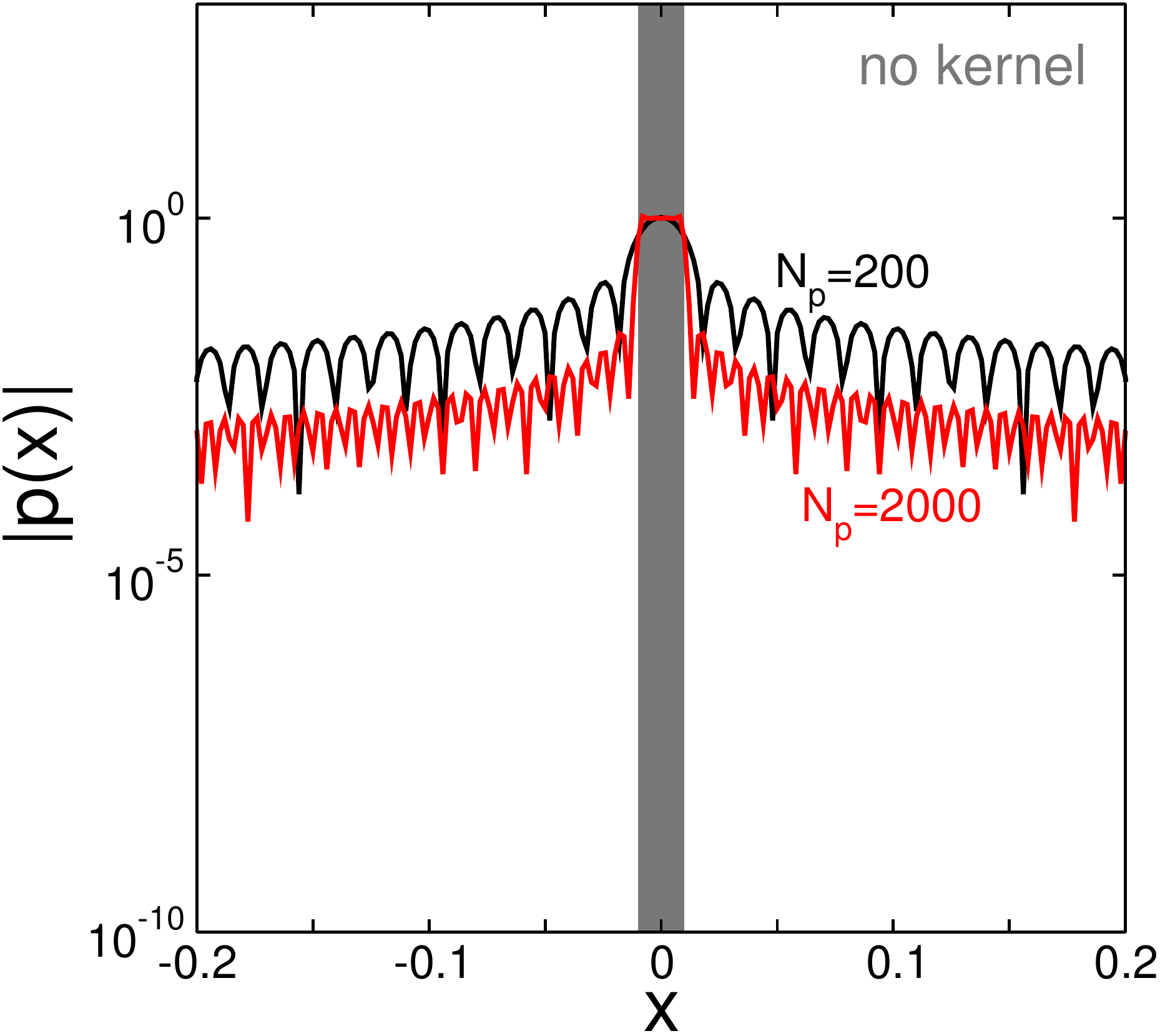}
\hspace*{\fill}
\includegraphics[width=0.4\textwidth]{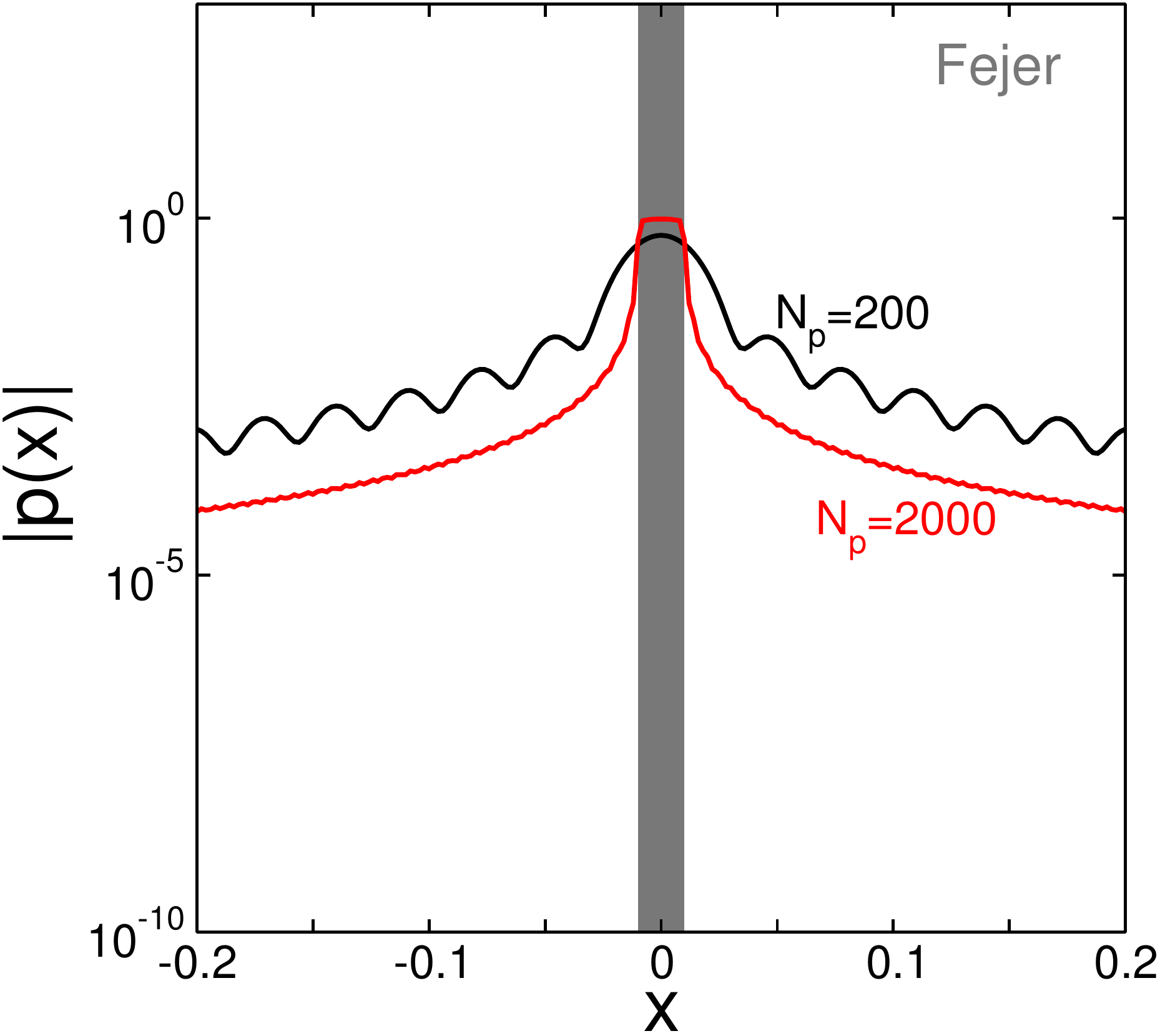}
\hspace*{\fill} \\[1ex]
\hspace*{\fill}
\includegraphics[width=0.4\textwidth]{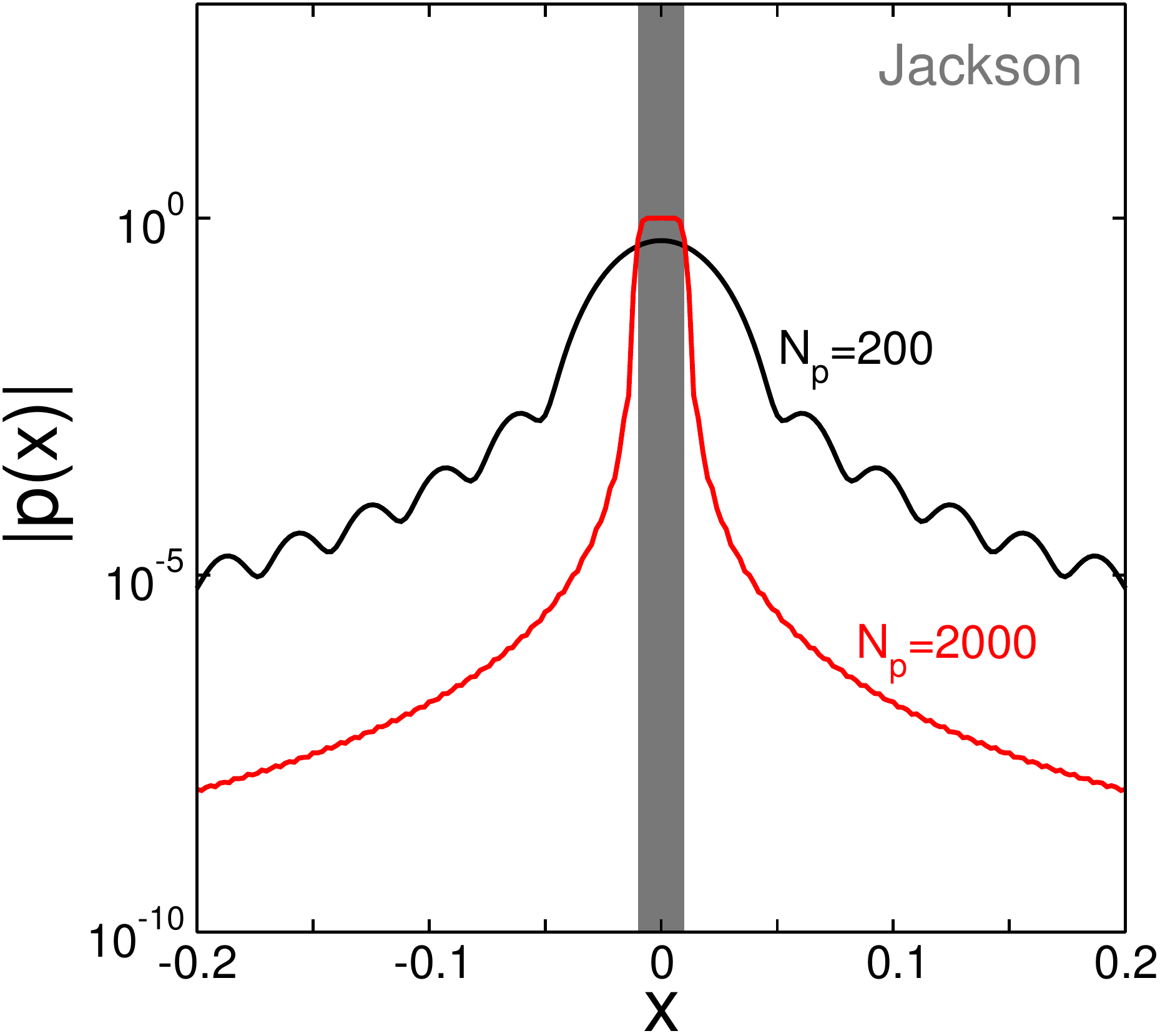}
\hspace*{\fill}
\includegraphics[width=0.4\textwidth]{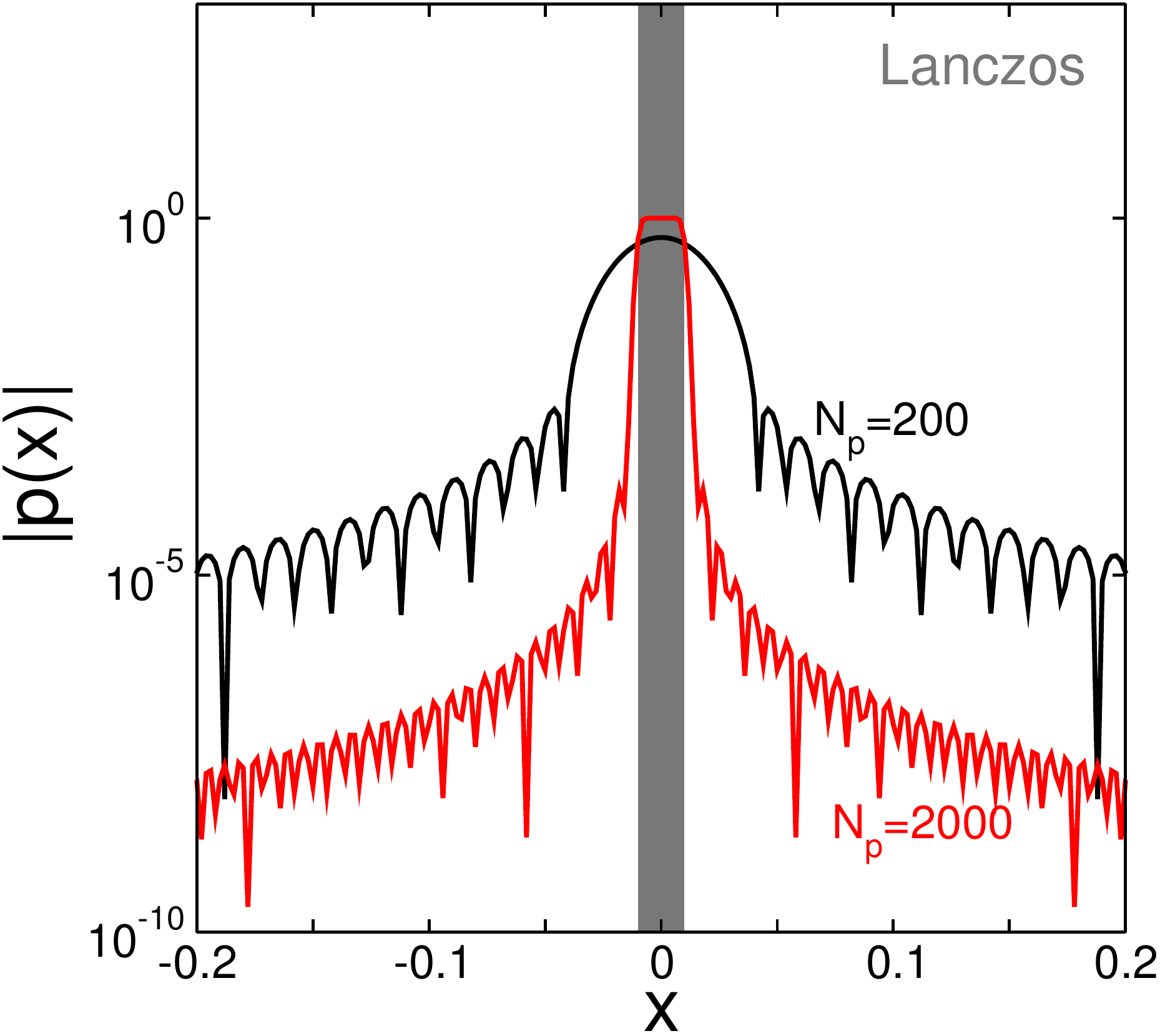}
\hspace*{\fill}
\caption{
Chebyshev polynomial approximations to the rectangular window function $W(x)$ for a target interval $[\lamu,\lamo]=[-0.01, 0.01] \subset [-1,1]=[a,b]$ (gray region), for polynomial degree $N_p=200$ (black curves) and $N_p=2000$ (red curves).
The approximations are shown without modification by a kernel ($g_n=1$ in Eq.~\eqref{W2}), and with the Fej\'{e}r, Jackson, and Lanczos ($\mu = 2$) kernels from Table~\ref{tab:Kernel}.
}
\label{fig:filters}
\end{figure}

Truncation of the infinite series~\eqref{W1} leads to a polynomial approximation to $W(x)$, which, however, suffers from Gibbs oscillations (see Fig.~\ref{fig:filters}).
The solution to the Gibbs phenomenon exploited in, e.g., KPM is to include suitable kernel damping factors $g_n$ in the expansion.
This gives the polynomial approximation
\begin{equation}\label{W2}
 W(x) \approx p_n(x) = \sum_{n=0}^{N_p} g_n c_n T_n( \alpha x + \beta)
\end{equation}
of degree $N_p$.

Three examples for the kernel coefficients $g_n$ are given in Table~\ref{tab:Kernel}, including the Jackson kernel~\cite{Jack12} known from KPM~\cite{WWAF06}.
The resulting approximations to $W(x)$ are plotted in Fig.~\ref{fig:filters}.
In ChebFD we will use the Lanczos kernel with $\mu=2$, which tends to give faster convergence 
than the other choices.

\begin{table}
\begin{tabular}{ccc}
  Fej\'{e}r & Jackson & Lanczos \\[0.5ex]\cline{1-3}\\[0.05ex]
  \multicolumn{1}{c|}{$\dfrac{N_p-n+1}{N_p+1}$} &  \multicolumn{1}{c|}{$ \dfrac1{N_p} \left( (N_p-n) \cos \dfrac{\pi n}{N_p} + \sin \dfrac{\pi n}{N_p} \cot \dfrac{\pi}{N_p} \right) $ }   & $\left( \dfrac{\sin (\pi n/(N_p+1))}{\pi n/(N_p+1)} \right)^\mu$
\end{tabular}
\caption{Coefficients $g_n$ of the Fej\'{e}r, Jackson, and Lanczos kernel (cf. Ref.~\cite{WWAF06}).
The Lanczos kernel has an additional parameter $\mu$ (later we will use $\mu=2$).}
\label{tab:Kernel}
\end{table}

Note that we have not yet defined what constitutes a `good'  filter polynomial $p(x)$.
Obviously, $|p(x)|$ should be (i) small for $x \not\in I_T$ and (ii) reasonably large for $x \in I_T$,
while (iii) keeping the degree of $p(x)$ small among all possible polynomials satisfying (i) and (ii).
We will examine these conditions in Sec.~\ref{sec:Theory} systematically.

Notice also that we construct filter polynomials starting from rectangular window functions, which is not necessarily the best approach. The fact that $W(x) = 1$ within the target interval is not relevant for convergence of ChebFD. Using window functions with $W(x) \ne 1$ that allow for better polynomial approximation can improve convergence~\cite{GKLAFPHKSWBRT16}. 
Another idea is the use of several different window functions, for example a comb of $\delta$-function peaks that cover the target interval.
Fortunately, as we will see now, the precise shape of the filter polynomials is not crucial for the overall convergence of ChebFD, although the convergence rate can be of course affected.

\begin{figure}
\hspace*{\fill}
\includegraphics[width=0.4\textwidth]{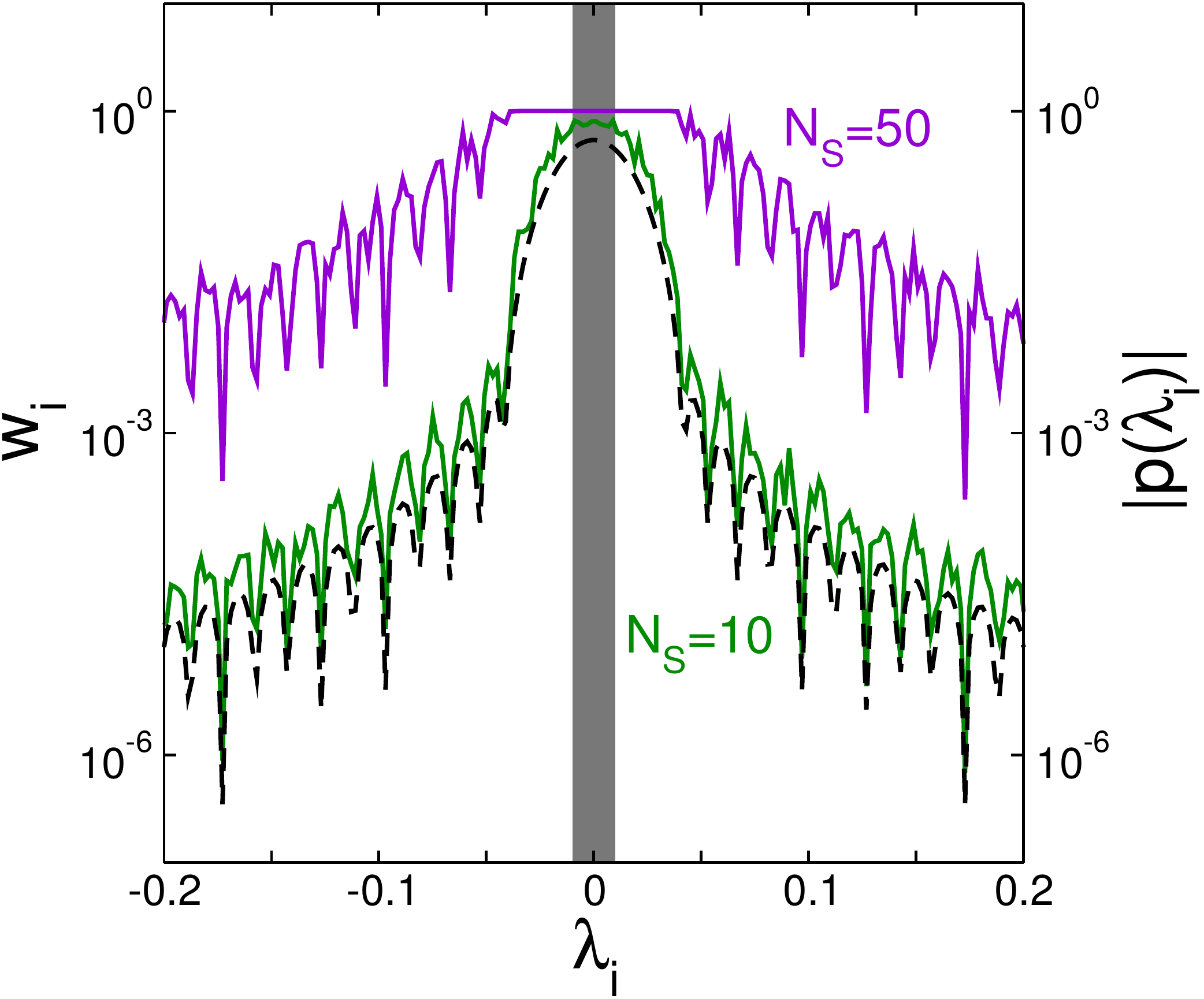}
\hspace*{\fill}
\includegraphics[width=0.4\textwidth]{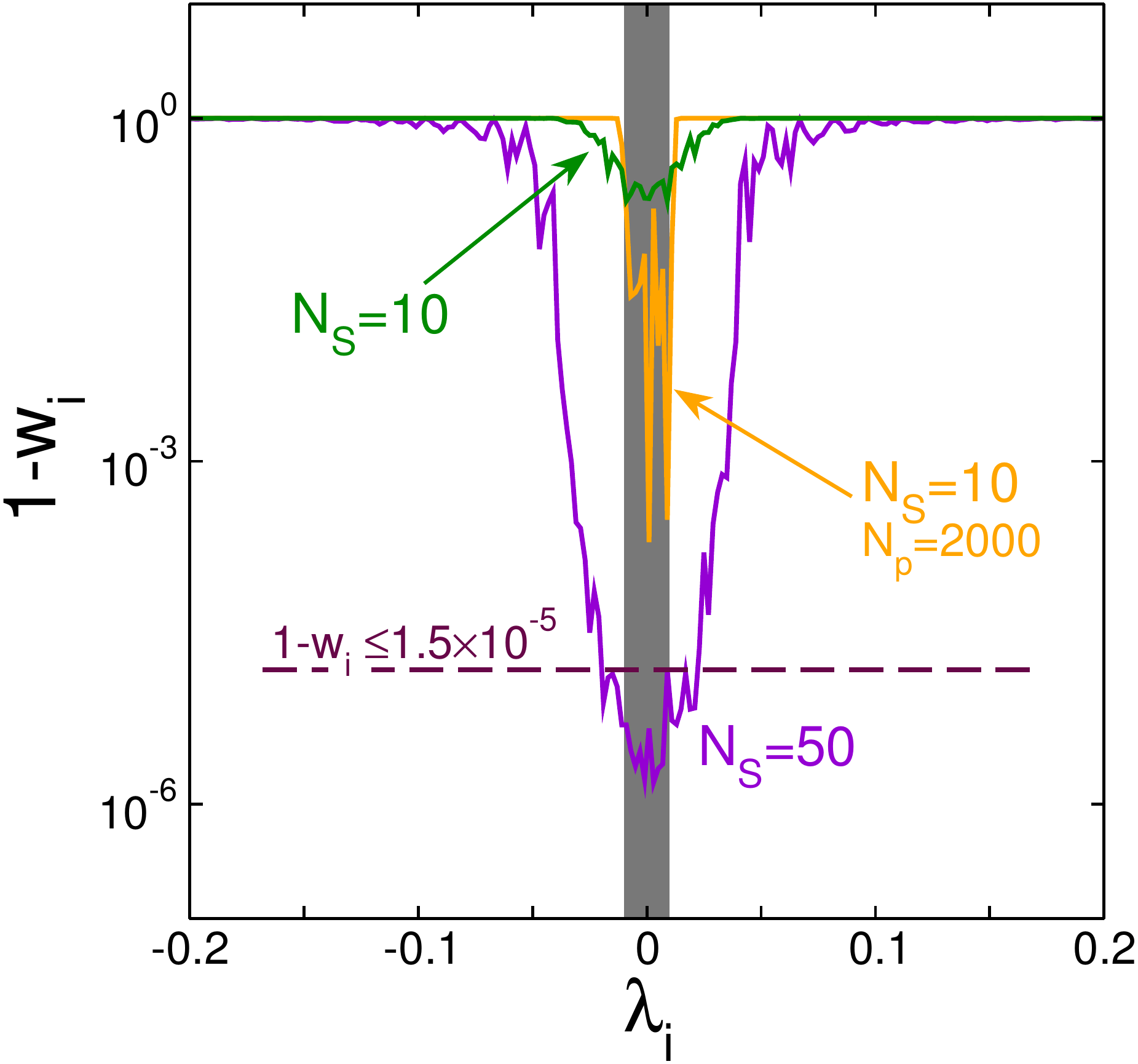}
\hspace*{\fill}
\caption{Weight $w_i$ of search vectors at eigenvectors (see Eq.~\eqref{weight1}) after one polynomial filtering step.
We use the filter polynomials for the target interval $[-0.01, 0.01] \subset [-1,1]$ (gray region) from Fig.~\ref{fig:filters} (here, for the Lanczos kernel with $\mu=2$).
The target space contains $N_T=10$ eigenvectors of the diagonal matrix described in the text, with dimension $D=10^3$. 
Left panel:
Weight $w_i$ for (i) $N_S=10$ (green curve) and (ii) $N_S=50$ (violet curve) search vectors, both for polynomial degree $N_p = 200$.
The black dashed curve shows the filter polynomial $|p(\lambda_i)|$.
Right panel:
Expected residual $1-w_i$ for cases (i), (ii) from the left panel,
and in addition for (iii) $N_S=10$ with $N_p=2000$ (orange curve).
The horizontal dashed line indicates the bound $1-w_i \le 1.5 \times 10^{-5}$ achieved within the target interval in case (ii).
}
\label{fig:filtermatrix}
\end{figure}

\subsection{Polynomial filtering: Numerical example}
\label{sec:Example1}

The effect of polynomial filtering is demonstrated in Fig.~\ref{fig:filtermatrix}.
We use the diagonal $D \times D$ matrix $H = \mathrm{diag}(\lambda_1, \dots, \lambda_D)$
with equidistant eigenvalues $\lambda_i = -1 + 2 i/(D+1)$.
The target interval $I_T = [-\delta,\delta]$ thus contains $N_T = \delta D$ eigenvectors of $H$.
We now apply the filter polynomial once,
and compare  cases with (i) a small number ($N_S=N_T$) and (ii) a large number ($N_S = 5 N_T \gg N_T$) of search vectors.
Shown is the weight 
\begin{equation}\label{weight1}
  w_i = \left( \sum_{k=1}^{N_S} |\langle \vec v_i , \vec y_k \rangle|^2 \right)^{1/2}
\end{equation}
of the filtered orthonormalized search vectors $\vec y_k$ at the eigenvectors $\vec v_i$ of $H$ (we use the Euclidean scalar product).
Notice that $w_i$ is invariant under orthogonal transformations of the subspace basis $\{ \vec y_1, \dots, \vec y_{N_S} \}$ of the search space formed by the orthonormalized vectors.
We have $0 \le w_i \le 1$, and $w_i=1$ if and only if the eigenvector $\vec v_i$ is contained in the filtered search space spanned by the $\vec y_k$.
If the search space contains an eigenvector approximation $\vec v$ to $\vec v_i$ with error $|\vec v - \vec v_i| \le \epsilon$, we have $w_i \ge 1- \epsilon$.
Conversely, $w_i \ge 1 - \epsilon$ implies that an eigenvector approximation with error below $2 \epsilon$ exists in the search space.
Values of $w_i$ close to one thus correspond to the existence of good eigenvector approximation among the filtered search vectors.
Finally, $w_i \ge 1-\epsilon$ for all $\lambda_i$ in the target interval indicates convergence to the desired accuracy $\epsilon$.

For a single search vector $\vec x = \sum_i x_i \vec v_i$ polynomial filtering gives
\begin{equation}\label{filterOne}
 p[H] \vec x = \sum_{i=1}^D x_i \, p(\lambda_i) \, \vec v_i \;,
\end{equation}
that is, $w_i \propto |p(\lambda_i)|$ on average.
The polynomial filter thus suppresses components outside of the target interval (where $p(\lambda_i) \approx 0$) in favor of components in the target interval (where $|p(\lambda_i)| \approx 1$).
This effect is clearly seen in the left panel of Fig.~\ref{fig:filtermatrix} for case (i) (with $N_S =N_T$).

The filter polynomial used in Fig.~\ref{fig:filtermatrix} is a rather `bad' approximation of the rectangular window function, despite the relatively large degree $N_p=200$.
One filter step with this polynomial suppresses components outside of the target interval only weakly.
Many iterations would be required in the ChebFD scheme, and convergence would be accordingly slow. 
One way to improve convergence is to improve the quality of the filter polynomial.
However, for narrow target intervals inside of the spectrum the degree of a `good' filter polynomial that leads to convergence in a few iterations had to become extremely large.

ChebFD follows a different approach:
Use a `moderately good' filter polynomial, even if it does not lead to fast convergence, and compensate by increasing the number of search vectors over the number of target vectors ($N_S \gg N_T$).
The rationale behind this approach is that polynomial filtering compresses the search space towards the target space, but the target space can accommodate only $N_T$ orthogonal vectors.
Therefore, orthogonalization of the filtered search vectors separates components inside and outside of the target space.
In this way ``overpopulating'' the target space with search vectors ($N_S \gg N_T$) can lead to faster convergence towards target vectors than Eq.~\eqref{filterOne} predicts (for $N_S \sim N_T$).

That this approach can significantly improve convergence is demonstrated by case (ii) (with $N_S = 5 N_T$) in Fig.~\ref{fig:filtermatrix}.
The same filter polynomial as in case (i) now achieves $w_i \approx 1$ throughout the target interval.
All target vectors are obtained from the search space with residual below $2 (1-w_i) \le 3 \times 10^{-5}$ in one filter step (see right panel).
The strength of this approach is also seen in comparison to case (iii) (right panel in Fig.~\ref{fig:filtermatrix}, with $N_S = N_T$ as in case (i), but $N_p = 2000$).
Although the filter polynomial is now much better than in cases (i), (ii) (see Fig.~\ref{fig:filters}, lower right panel) the residual weight $1-w_i$ has not diminished substantially.
The accuracy achieved in case (iii) is still far from that achieved in case (ii), 
although the number of spMVMs has doubled ($N_S \times N_p = 10 \times 2000$ versus $50 \times 200$). 
Evidently, there is a trade-off between search space size and filter polynomial degree
that requires further analysis.

\subsection{Theoretical analysis of filter polynomial quality and ChebFD convergence}
\label{sec:Theory}

The previous numerical experiment illustrates the essential mechanism behind polynomial filtering in ChebFD.
Our theoretical analysis of the ChebFD convergence is based on a measure for the quality of the filter polynomial that depends on the target interval and search interval size.

Arrange the eigenvalues $\lambda_i$ of $H$ such that $|p(\lambda_1)| \ge |p(\lambda_2)| \ge \dots \ge |p(\lambda_D)|$.
By construction, the filter polynomial is largest inside of the target interval $I_T$, such that the target eigenvalues appear as $\lambda_1, \dots, \lambda_{N_T}$.
Therefore, through polynomial filtering the components of the search vectors that lie in the target space are multiplied by a factor $\ge |p(\lambda_{N_T})|$.
On the other hand, because we use $N_S > N_T$ search vectors, components perpendicular to the target space are multiplied by a factor $\le p(\lambda_{N_S})$.
The damping factor of the filter polynomial, by which the unwanted vector components in the search space are reduced in each iteration, is $\sigma = |p(\lambda_{N_S})/p(\lambda_{N_T})|$.

As seen in Fig.~\ref{fig:filters} the filter polynomials decay with distance from the target interval.
We will, therefore, assume that the eigenvalues $\lambda_{N_T+1}, \dots, \lambda_{N_S}$ in the above list are those closest to the target interval.
With this assumption we arrive at the expression
\begin{equation}
 \sigma = \frac{\max\limits_{x \in [a,b], x \not \in I_S} |p(x)|}{\min\limits_{x \in I_T} |p(x)| }
\end{equation}
for the damping factor of the filter polynomial,
where $I_S = [\underline\nu,\overline\nu] \supset I_T = [\lamu,\lamo]$ is the search interval
that contains the eigenvalues $\lambda_1, \dots, \lambda_{N_S}$.
For a symmetrically shaped polynomial filter 
the search interval evenly extends the target interval 
by a margin $\delta' > 0$,
i.e., $\underline\nu = \lamu - \delta'$, $\overline\nu = \lamo + \delta'$  (see Fig.~\ref{fig:sketch}).

If the filter polynomial $p(x)$ is applied $M$ times (i.e., the filter polynomial used is $p(x)^M$),
the damping factor is $\sigma^M$ and the total polynomial degree is $M \times N_p$.
To achieve a reduction of the unwanted vector components outside of the target space by a factor $\epsilon$, therefore, requires about $N_\mathrm{iter} \approx \log_{10} \epsilon  / \log_{10} \sigma$ iterations in the ChebFD scheme, 
or $\eta \times (-\log_{10} \epsilon)$ spMVMs for each vector, where
\begin{equation}
 \eta = - \frac{N_p}{\log_{10} \sigma}
\end{equation}
is our (lower-is-better) measure for the quality of the filter polynomial\footnote{Note that the number of iterations $N_\mathrm{iter}$ is an integer, therefore the actual number of spMVMs required in practice will be somewhat larger than this estimate. Such details do not change the general analysis.}.

\begin{figure}
\hspace*{\fill}
\includegraphics[width=0.4\textwidth]{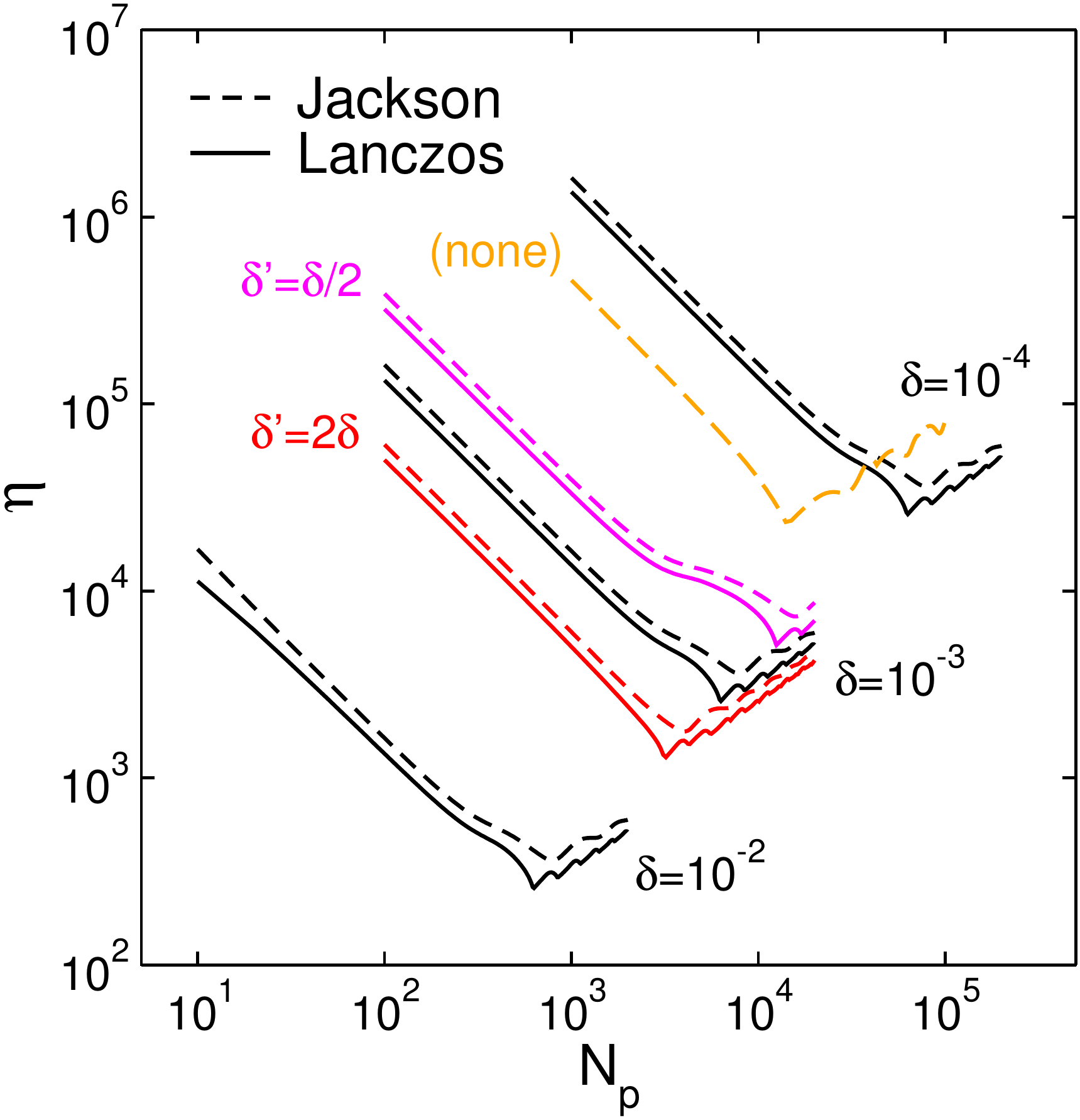}
\hspace*{\fill}
\includegraphics[width=0.4\textwidth]{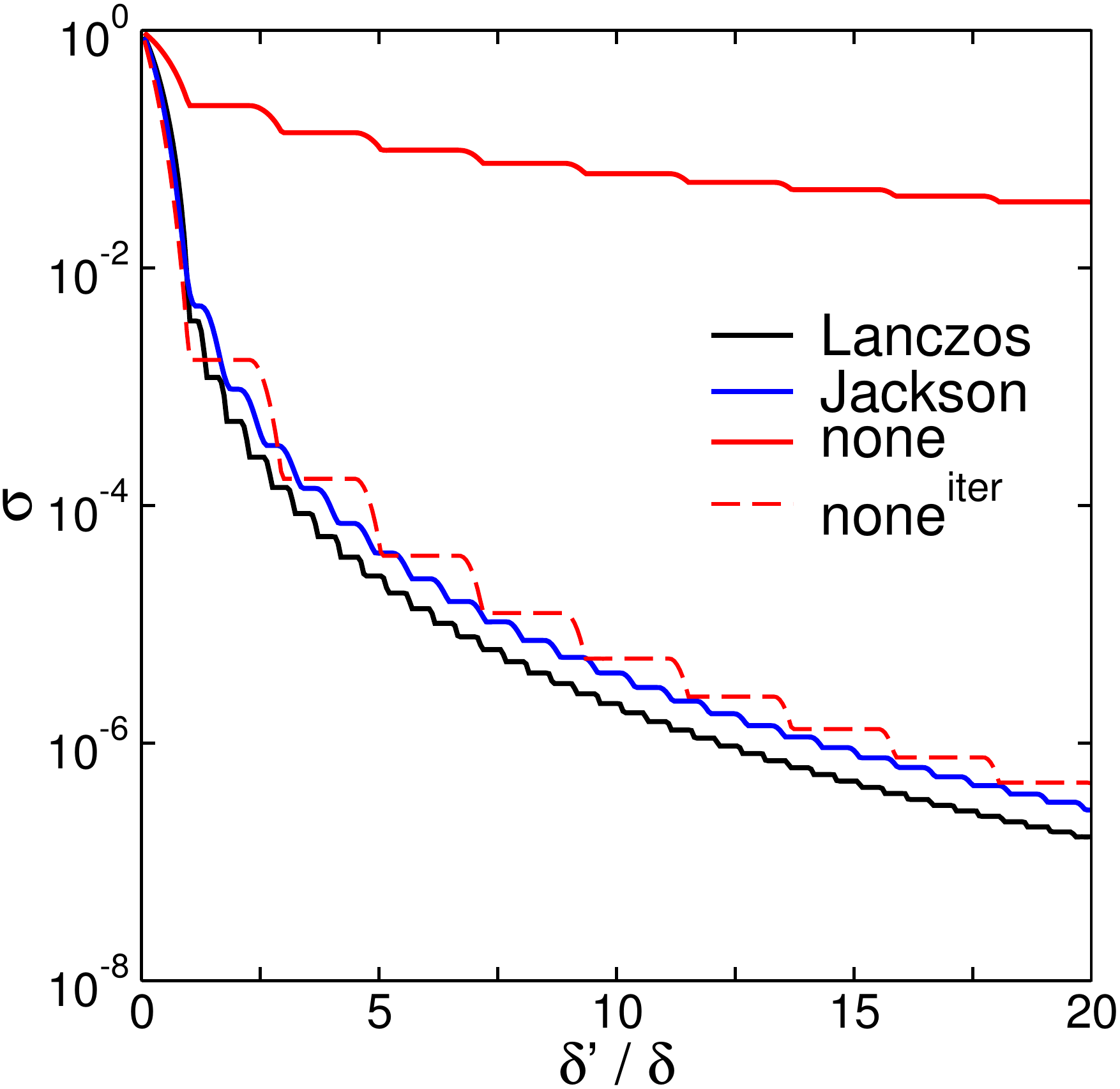}
\hspace*{\fill}
\caption{Left panel: Filter polynomial quality $\eta$ as a function of $N_p$,
for target interval half width $\delta = 10^{-2}, 10^{-3}, 10^{-4}$ and search interval margin $\delta' = \delta$,
with the Jackson and Lanczos ($\mu=2$) kernel.
For $\delta = 10^{-3}$, additional curves with $\delta ' = \delta/2, 2 \delta$ are shown,
and for $\delta = 10^{-4}$ the additional curve marked ``(none)'' gives the result obtained without kernel factors.
Right panel: Damping factor $\sigma$ as a function of the search interval margin $\delta'$, for $\delta = 10^{-3}$ and the optimal polynomial degree $N_p$  as determined in the left panel for $\delta'=\delta$. It is $N_p =6251, 7899, 1424$ for the Lanczos, Jackson, no kernel.
The curve ``none (iter)'' gives the result for the $6251/1424 \approx 4.4$-th power of the ``none'' filter polynomial.
}
\label{fig:polyeff1}
\end{figure}

The behavior of the filter polynomial quality $\eta$ is plotted in Figs.~\ref{fig:polyeff1},~\ref{fig:polyeff2}.
We concentrate here, and in the following, on the symmetric situation
$I_T = [-\delta,\delta]$ and $I_S = [-(\delta+\delta'),\delta+\delta']$ in the center of the interval $[-1,1]$.
Similar considerations apply, with modified numerical values, also to situations with $\lamu \ne - \lamo$ as long as the target interval is located inside of the spectrum.
For the computation of extreme eigenvalues other considerations would apply.

For all combinations of $\delta$, $\delta'$ there is an optimal value of $N_p$ that leads to minimal $\eta$.
As seen in Fig.~\ref{fig:polyeff1} (left panel) the optimal $N_p$ and $\eta$ scale as $\delta^{-1}$ for fixed $\delta'/\delta$, and as $\delta'^{-1}$ for fixed $\delta$.
We observe that high polynomial degrees $N_p$ are required already for relatively large values of $\delta$.
The Lanczos kernel (with $\mu=2$) performs always better than the Jackson kernel, for which the broadening of the window function is too pronounced (cf. Fig.~\ref{fig:filters}).
In some situations the filter polynomial without kernel modifications (``none'' kernel), for which the optimal $\eta$ occurs at significantly smaller $N_p$ because the width of the window function approximation remains small,  is even better.

As seen in Fig.~\ref{fig:polyeff1} (right panel) the kernels differ in their decay away from the target interval (see also Fig.~\ref{fig:filters}).
Even if the ``none'' kernel can outperform the Lanczos kernel in terms of $\eta$, it does not decay as fast with distance from the target interval.
In the application, where we do not know the search interval precisely,
the robustness of the filter polynomial under changes of the search interval becomes important. Because of the faster and smoother decay---and for this reason only---we will use the Lanczos kernel ($\mu=2$) in the following and not the ``none'' kernel.
In any case, the construction of better filter polynomials remains worthy of investigation.

\subsection{Optimal choice of algorithmic parameters}
\label{sec:Optimal}

Based on the theoretical analysis from the previous subsection we will now establish practical rules for the choice of the two algorithmic parameters $N_p$ \textit{(filter polynomial degree)} and $N_S$ \textit{(number of search vectors)}. 

\begin{figure}
\hspace*{\fill}
\includegraphics[width=0.4\textwidth]{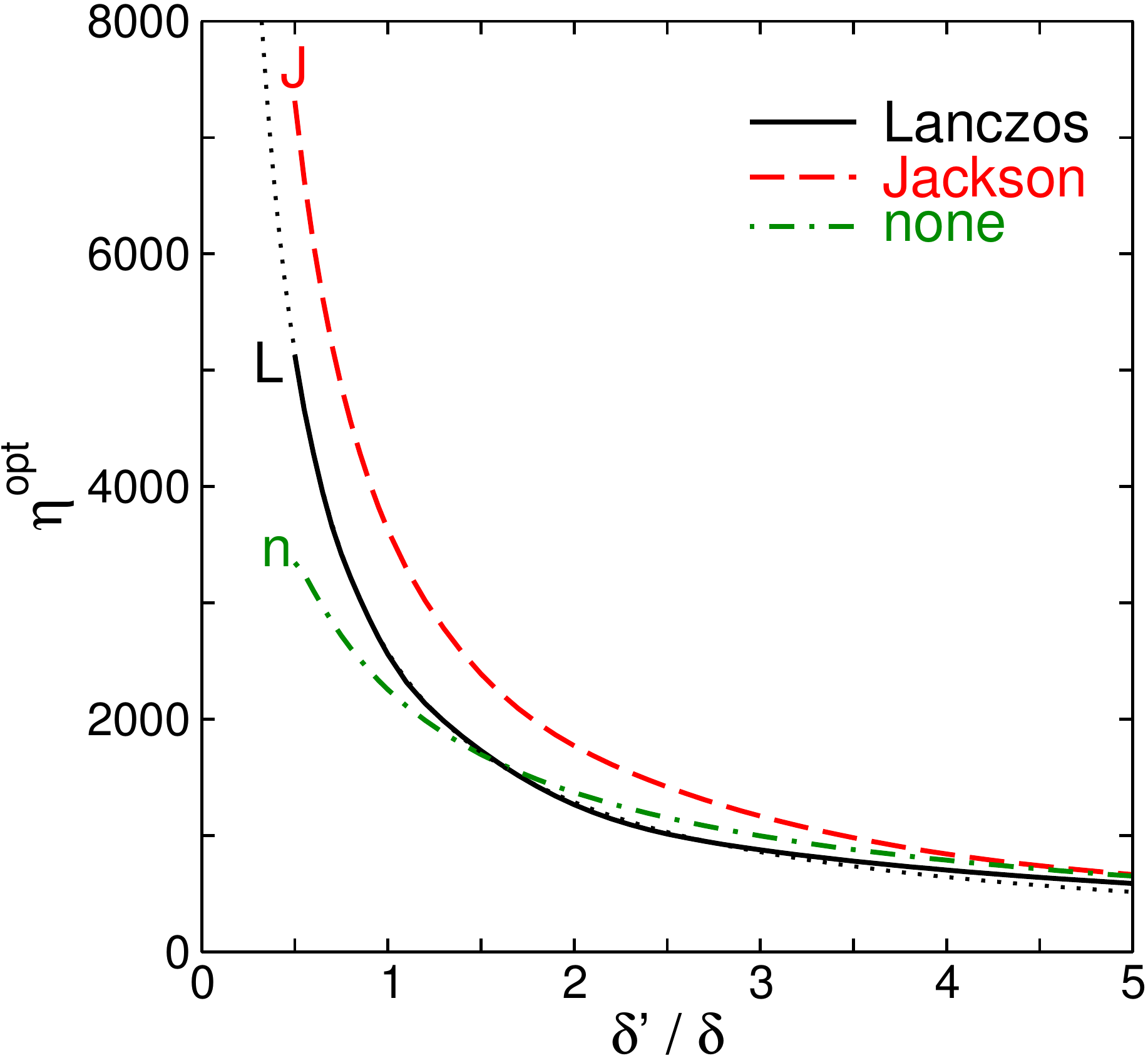}
\hspace*{\fill}
\includegraphics[width=0.4\textwidth]{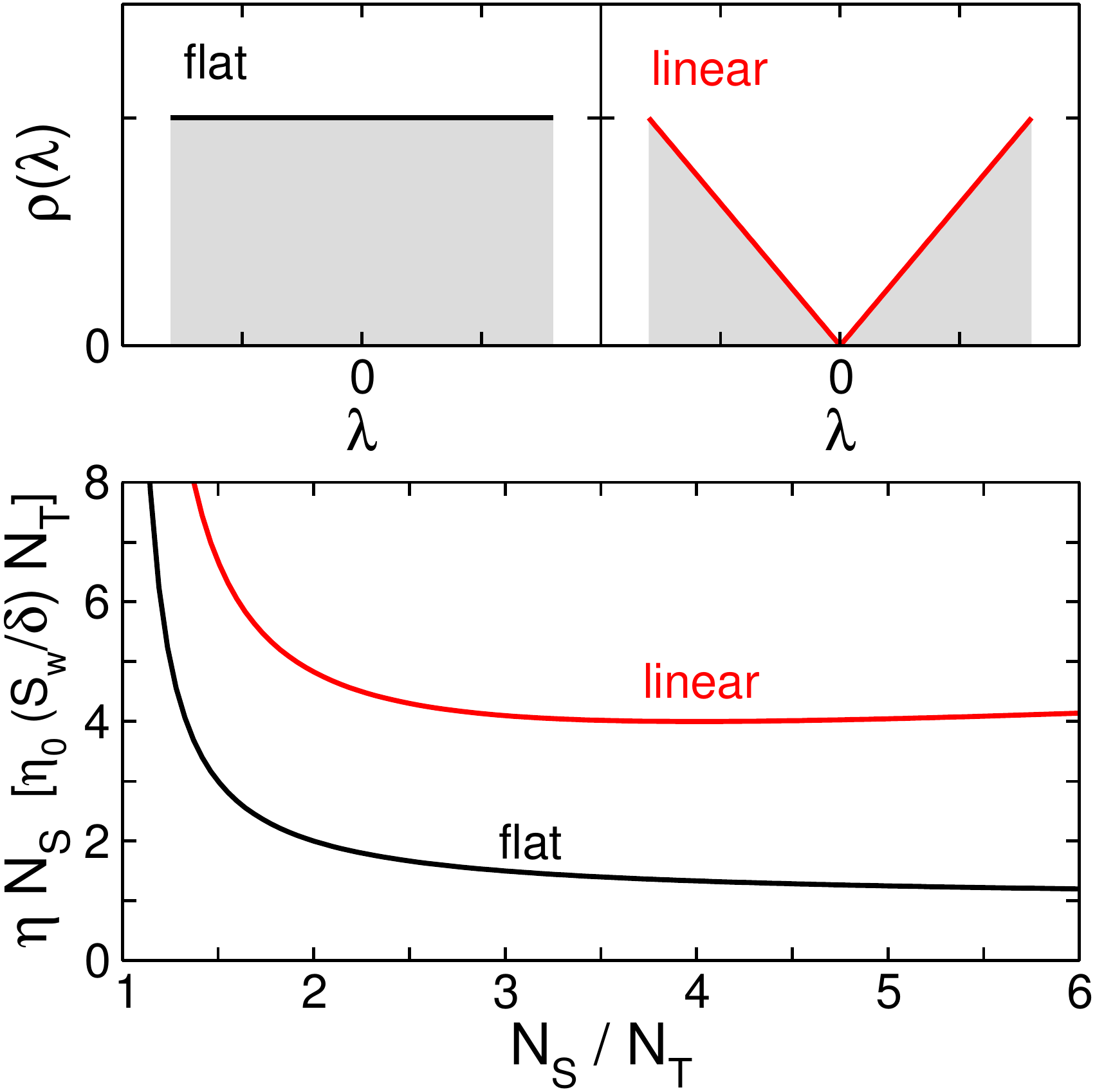}
\hspace*{\fill}
\caption{Left panel: Value of the optimal (i.e., minimal possible) filter polynomial quality $\eta^\mathrm{opt}$ as a function of search interval margin $\delta'/\delta$, for $\delta = 10^{-3}$.
The dotted curve gives the fit $\eta \propto 1/\delta'$ from Eq.~\eqref{etafit}.
Right panel: Estimate for the effort $\eta N_S$, measured in units of $\eta_0 (S_w/\delta) N_T$, for a flat and linear density of states according to Eqs.~\eqref{etanflat},~\eqref{etanlinear}.
}
\label{fig:polyeff2}
\end{figure}

Because the filter polynomial is applied to all search vectors,
the total number of spMVMs is $N_\mathrm{iter} \times N_p \times N_S$,
or $\eta \times N_S \times (-\log_{10} \epsilon)$.
The last factor in this expression is specified by the user,
so that we will try to minimize the product $\eta \times N_S$.
The central information used here is the scaling
\begin{equation}\label{etafit}
 \eta^\mathrm{opt} \simeq \eta_0  \frac{S_w}{\delta} \frac{\delta}{\delta'}
\end{equation}
of the optimal filter polynomial quality (cf. Fig.~\ref{fig:polyeff2})
with the half width $S_w=(b-a)/2$ of the interval $[a,b]$ that contains the spectrum, the target interval half width $\delta$, and the search interval margin $\delta'$.
The optimal polynomial degree scales as  $N_p^\mathrm{opt} \simeq N_0 (S_w/\delta) (\delta/\delta')$
(cf. Fig.~\ref{fig:polyeff1}).
Only the dependence on $\delta'$ is relevant for the following argument.
We introduce the additional factor ``$\delta / \delta$'', which cancels in Eq.~\eqref{etafit},
to obtain the product of the relative width $S_w / \delta$ (of the target interval) and $\delta/\delta'$ (of the search interval margin).

The approximate scaling relations for $\eta^\mathrm{opt}$ and $N_p^\mathrm{opt}$ express the fact that the width\footnote{\label{afootnote}For example, the variance of the Jackson kernel is $\pi/ N_p$ in the center of the interval $[-1,1]$ (see Ref.~\cite{WWAF06}), which gives the full width at half-maximum $\sqrt{2 \ln  2} \, \pi (b-a)/N_p \approx 3.7 (b-a)/N_p$ of the Gaussian function that is obtained from a $\delta$-peak through convolution with the kernel.} of the (Lanczos or Jackson) kernel scales as $S_w / N_p$.
Therefore, if the filter polynomial from Eq.~\eqref{W2} should decrease from $|p(x)| \approx 1$ inside of the target interval to $|p(x)| \ll 1$ outside of the search interval, i.e., over the distance $\delta'$, the polynomial degree has to scale as $S_w / \delta'$.
The above scaling relations follow.
The constants $\eta_0$, $N_0$ in these relations depend on the kernel, and can be obtained from a fit as in Fig.~\ref{fig:polyeff2}. For the Lanczos kernel (with $\mu=2)$ we have $\eta_0 = 2.58$, $N_0 = 6.23$ near the center of the spectrum (cf. Table~\ref{tab:lanczoskernel}).

The optimal choice of $N_S$ now depends on the distribution of eigenvalues outside of the target interval,
which can be expressed by the eigenvalue count, or density of states (DOS), $\rho(\lambda) = \sum_i \delta(\lambda -\lambda_i)$.
By definition, $D = \int_{a}^b \rho(\lambda) d\lambda$ gives the total number of states, i.e., the matrix dimension.
The number of target vectors is given by
\begin{equation}\label{DOS:target}
N_T = \int_{\lamu}^{\lamo} \rho(\lambda) d\lambda \;,
\end{equation}
while the search interval $I_S = [\underline\nu,\overline\nu]$ can be determined from the condition
\begin{equation}\label{DOS:search}
 N_S = \int\limits_{\underline\nu}^{\overline\nu} \rho(\lambda) \mathrm d\lambda \;,
\end{equation}
where we can assume that $\underline\nu = \lamu - \delta'$, $\overline\nu = \lamo + \delta'$.

Two situations are of particular practical importance:
A flat density of states, which stays constant in the vicinity of the target interval, and a linearly increasing density of states with a `pseudo'-gap in the target interval.
The first situation occurs, e.g., in Anderson localization,
the second situation occurs for the graphene and topological insulator materials.
The optimal choice in these cases, as depicted in Fig.~\ref{fig:polyeff2}, is as follows.

\paragraph{Flat DOS}
For a flat DOS $\rho(\lambda) = \rho_0 \equiv \text{const.}$
we have $N_T = 2 \rho_0 \delta$ and $N_S = 2 \rho_0 (\delta + \delta')$.
Then, $\delta' / \delta = N_S/N_T - 1$ or
\begin{equation}\label{etanflat}
\eta N_S = \eta_0 \frac{S_w}{\delta}   N_T  \big(1 - N_T/N_S \big)^{-1} \;.
\end{equation}
This function decreases monotonically with $N_S$,
with $\eta N_S \to \eta_0 (S_w/\delta) N_T$ for $N_S \gg N_T$. 

\paragraph{Linear DOS}
For a linearly increasing DOS $\rho(\lambda) = \rho_1 |\lambda|$
we have $N_T = \rho_1 \delta^2$ and $N_S = \rho_1 (\delta + \delta')^2$.
Now, $\delta '/ \delta = \sqrt{N_S/N_T} - 1$ or
\begin{equation}\label{etanlinear}
\eta N_S =\eta_0 \frac{S_w}{\delta} N_T  \frac{\sqrt{N_S/N_T}}{1-\sqrt{N_T/N_S}} \;.
\end{equation}
This function is minimal at $N_S = 4 N_T$, where $\eta N_S = 4 \eta_0 (S_w/\delta) N_T$.

\medskip

It is instructive to compare the scaling of the effort with the matrix dimension $D$ in both cases.
For the flat DOS, with $D \simeq 2 \rho_0 S_w$, the scaling is
$\eta N_S \simeq \eta_0 D$.
For the linear DOS, with $D \simeq \rho_1 S_w^2$, the scaling is
$\eta N_S \simeq 4 \eta_0 \sqrt{D N_T}$.
For the linear DOS the required total effort is smaller by a factor $4 \sqrt{N_T/D}$.

Note that our considerations apply under the assumption that $N_S / D$ is a small but non-negligible ratio, such that the DOS can be approximated by a continuous function, and that the target interval lies inside of the spectrum.
They break down for $N_S \to 1$, where the discreteness of eigenvalues comes into play,
and at the edges of the spectrum where other scaling relations hold.

\section{The Chebyshev filter diagonalization scheme}
\label{sec:ChebFD}

\begin{table}
\centering
\small
\begin{tabular}{*{11}{c}}
& \multicolumn{10}{c}{interval center $c / S_w$}\\
& 0 & 0.1 & 0.2 & 0.3 & 0.4 & 0.5 & 0.6 & 0.7 & 0.8 & 0.9 \\\cline{2-11}
$N_0$ \rule{0pt}{3.2ex} &  6.23 & 6.20 & 6.10 & 5.94 & 5.71 &  5.40 & 4.99 & 4.46 & 3.75 & 2.73 \\
$\eta_0$ & 2.58 & 2.57 & 2.53 & 2.46 & 2.37 & 2.24 & 2.07 & 1.85 & 1.55 & 1.13 
\end{tabular}
\caption{Values for $\eta_0$ and $N_0$ in the approximate scaling relations
(cf. Eq.~\eqref{etafit})  for the Lanczos kernel ($\mu=2$), for various positions $c$ of the target interval $[c-\delta,c+\delta]$ in the spectral interval $[-S_w, S_w]$. The values reported have been determined for $\delta/S_w =10^{-3}$. Notice that the scaling relations hold for a large range of $\delta$-$\delta'$ combinations, although only approximately.
Notice further that $N_0$ and $\eta_0$ decrease towards the edges of the spectrum ($c \to \pm S_w$),
that is, the minimal required effort becomes smaller.}
\label{tab:lanczoskernel}
\end{table}

The ChebFD scheme is the straightforward realization of the idea sketched in the previous section:
(Over-) populate the target space with many filtered search vectors.
Fast convergence is not achieved by unlimited improvement of the filter polynomial quality
but relies essentially on the use of many search vectors.
The requirements for memory resources are accordingly high, which suggests that ChebFD is best deployed in a large-scale computing environment.

\subsection{Algorithm}

The ChebFD scheme comprises the following steps:

\begin{enumerate}

\item Determine parameters $a,b$ such that the spectrum of $H$ is contained in the interval $[a,b]$. 
Estimate the DOS $\rho(\lambda)$.

\item Estimate the number of target vectors $N_T$ in the specified target interval. Choose the number of search vectors $N_S \gtrsim 2 N_T$ according to this estimate.

\item Estimate the size of the search interval from $N_S$. Choose the filter polynomial degree $N_p$ according to this estimate.

\item Construct $N_S$ random search vectors $\vec x_1, \dots, \vec x_{N_S}$. 

\rule{\linewidth}{0.5pt}

\item Apply the polynomial filter: $\vec y_k = p[H] \vec x_k$ (see Fig.~\ref{fig:ChebFD}).
\label{corestep}

\item Orthogonalize the filtered search vectors $\{ \vec y_k \}$.

\item Compute the Rayleigh-Ritz projection matrix $\langle \vec y_k,  H \vec y_l \rangle $ and the corresponding Ritz pairs $(\tilde \lambda_k, \tilde {\vec v}_k)$, together with their residuals $\vec r_k$.

\item Check for convergence: 

\begin{enumerate}
\item Exit if $|\vec r_k| \le \epsilon$ for all Ritz values $\tilde \lambda_k$ in the target interval.
\item Otherwise: Restart at step~\ref{corestep} with the new search vectors $\{ \vec y_k \}$.
\end{enumerate}

\end{enumerate}

\begin{figure}
Step~\ref{corestep}: Apply the polynomial filter \\[-1ex]
\rule{0.5\textwidth}{0.3mm}
\begin{algorithmic}[1]  
\For{ $k=1$\ to\ $N_S$ } \Comment{{\it First two recurrence steps}}
\State $\vec u_k = (\alpha H + \beta \mathbbm 1) \vec x_k$ \Comment{\texttt{spmv()}}
\State $\vec w_k = 2 (\alpha H + \beta \mathbbm 1) \vec u_k - \vec x_k$ \Comment{\texttt{spmv()}}
\State $\vec x_k = g_0 c_0 \vec x_k + g_1 c_1 \vec u_k + g_2 c_2 \vec w_k$ \Comment{\texttt{axpy \& scal}}
\EndFor
\For{ $n=3$\ to\ $N_p$ } \Comment{{\it Remaining recurrence steps}}
\For{ $k=1$\ to\ $N_S$ }
\State $\mbox{swap} ( \vec w_k ,  \vec u_k )$ \Comment{swap pointers}
\State $\vec w_k = 2 (\alpha H + \beta \mathbbm 1) \vec u_k -  \vec w_k $ 
\Comment{\texttt{spmv()}}
\State $\vec x_k = \vec x_k + g_n c_n \vec w_k$ \Comment{\texttt{axpy}}
\EndFor
\EndFor
\end{algorithmic}
\caption{The computational core (step~\ref{corestep}) of the ChebFD scheme: replace $\vec x_k$ by $p[H] \vec x_k$.
The two loops over $k$ can each be replaced by a single operation on block vectors of size $n_b = N_S$ (see Sec.~\ref{sec:Impl}).
}
\label{fig:ChebFD}
\end{figure}

Steps 1--4 comprise the preparatory phase, and steps 5--8 the iterative core phase of the ChebFD scheme.
A few remarks on the individual steps of the algorithm:

\paragraph{Step 1}
The interval $[a,b]$ can be computed with a few ($20-30$) Lanczos iterations,
which we prefer over the crude estimates from, e.g., Gershgorin's theorem.
Good estimates for the DOS can be computed with the KPM in combination with stochastic sampling of the matrix trace (see Ref.~\cite{WWAF06} for details).
This strategy has been known to physicists for a long time~\cite{SRVK96}  
and has recently been adopted by researchers from other communities~\cite{dNPS14}.
 
 \paragraph{Step 2}
The number of target vectors can be estimated from the DOS
according to Eq.~\eqref{DOS:target}.
Following the analysis from Secs.~\ref{sec:Theory},~\ref{sec:Optimal}
the number of search vectors should be a small multiple of the number of target vectors.
We recommend using $2 \le N_S/N_T \le 4$.
 
 \paragraph{Step 3}
 
The search interval width can be estimated from the DOS
according to Eq.~\eqref{DOS:search}.
The optimal value of $N_p$ follows from minimization of $\eta$ for the given interval configuration,
as in Sec.~\ref{sec:Theory}.
The minimization does not require the matrix $H$, can be performed prior to the ChebFD iteration, and requires negligible additional effort.
For the Lanczos kernel the choice of $N_p$ can be based on the values in Table~\ref{tab:lanczoskernel}.
With the corresponding value of $\eta$ the expected computational effort can be deduced already at this point.

\paragraph{Step 5}
The computational core of ChebFD is the application of the filter polynomial to the search space in step~\ref{corestep},
as depicted in Fig.~\ref{fig:ChebFD}.
Only during this step in the inner iteration cycle (steps 5--8) the matrix $H$ is addressed through spMVMs.
Therefore, our performance engineering efforts (Sec.~\ref{sec:Impl}) focus on this step.

\paragraph{Step 6}
Orthogonalization should be performed with a rank-revealing technique such as SVQB~\cite{SW02} or TSQR~\cite{DGHL12}.  
Because the damping factors $\sigma$ of the filter polynomials used in practice are still large compared to machine precision, the condition number of the Gram matrix of the filtered search vectors 
usually remains small. Therefore we can use SVQB, which is not only simpler to implement than TSQR but also more easily adapted to the row-major storage of vectors used in our spMVM (see Sec.~\ref{sec:Impl}).
For the computation of the required scalar products we use Kahan summation~\cite{Kah65}
 to preserve accuracy even for long vectors (i.e., large $D$).
If orthogonalization reduces the size of the search space we can add new random vectors to keep $N_S$ constant.

\paragraph{Step 7}
The Ritz pairs computed in step 7 will contain a few `ghost' Ritz pairs, for which the Ritz values lie inside of the target interval but the residual is large and does not decrease during iteration.
To separate these `ghosts' from the Ritz vectors that do converge to target vectors we discard all Ritz pairs for which the residual lies above a certain threshold (here: if it remains larger than $\sqrt{\epsilon}$).
A theoretically cleaner approach would be the use of harmonic Ritz values~\cite{Mor91,PPV95},
but they are not easily computed within the present scheme, and we did not find that they speed up convergence or improve robustness of the algorithm.
Therefore, we decided to use ordinary Ritz values together with the above acceptance criterion in the present large-scale computations.

\paragraph{Step 8}
The convergence criterion stated here uses the parameter $\epsilon$ directly as specified by the user, similar to the criteria used in other standard eigensolvers~\cite{ARPACK}.
Perturbation theory explains how $\epsilon$ is related to the error of the eigenvalues and eigenvectors~\cite{BhatiaPerturbation}.
For eigenvalues, a small residual $|\vec r_k | \le \epsilon$ guarantees a small error
$|\lambda_k - \tilde\lambda_k| \le \epsilon \| H \|$ of the corresponding eigenvalue
(the norm $\|H \|$ is known from step 1).
If precise eigenvectors are required the convergence check should also take the distance to neighboring eigenvalues into account, and can be performed by a user-supplied routine that accesses all the Ritz pairs and residuals.

\subsection{Parameter selection: Numerical experiment}

\begin{figure}
\hspace*{\fill}
\includegraphics[width=0.4\textwidth]{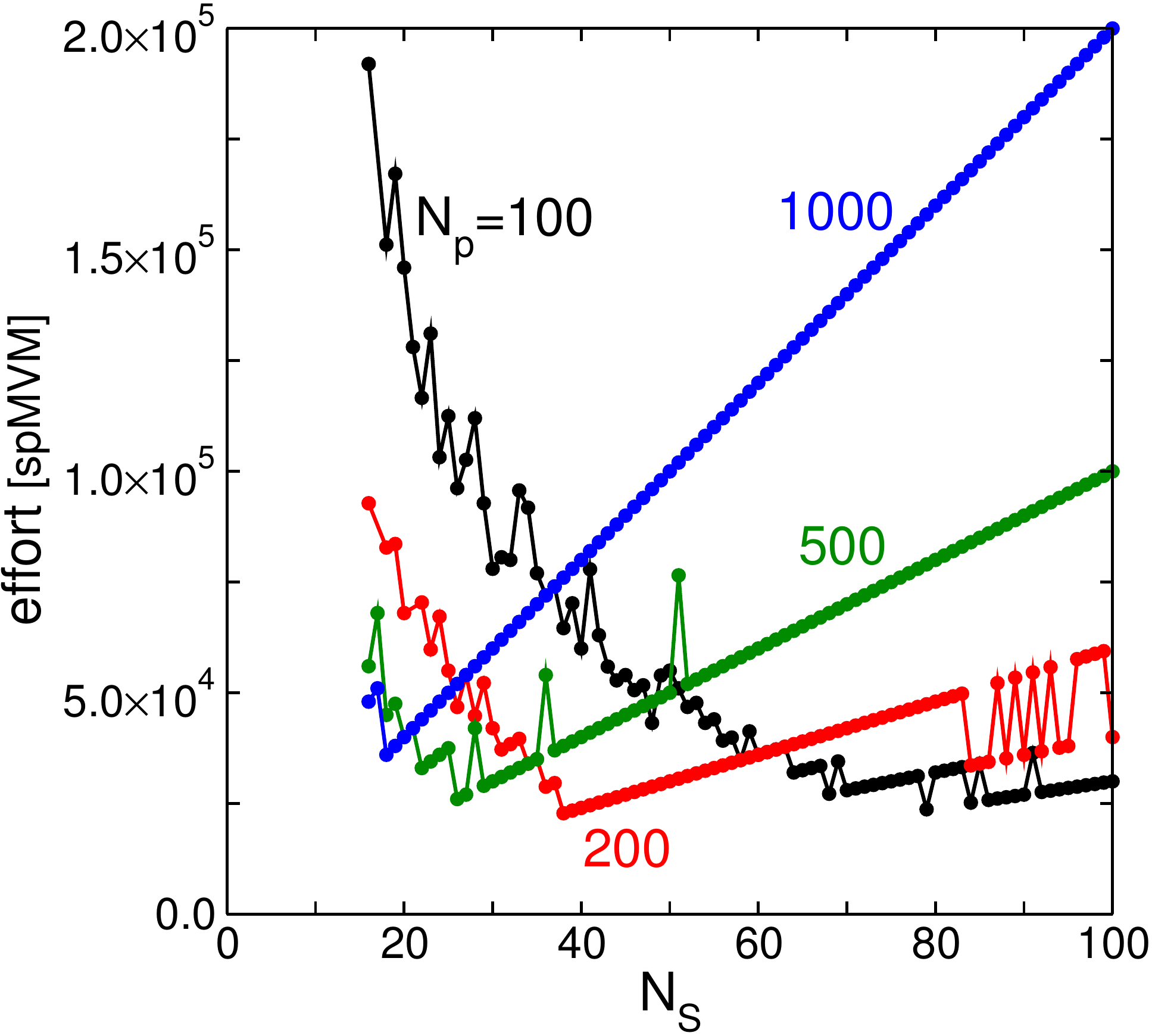}
\hspace*{\fill}
\includegraphics[width=0.4\textwidth]{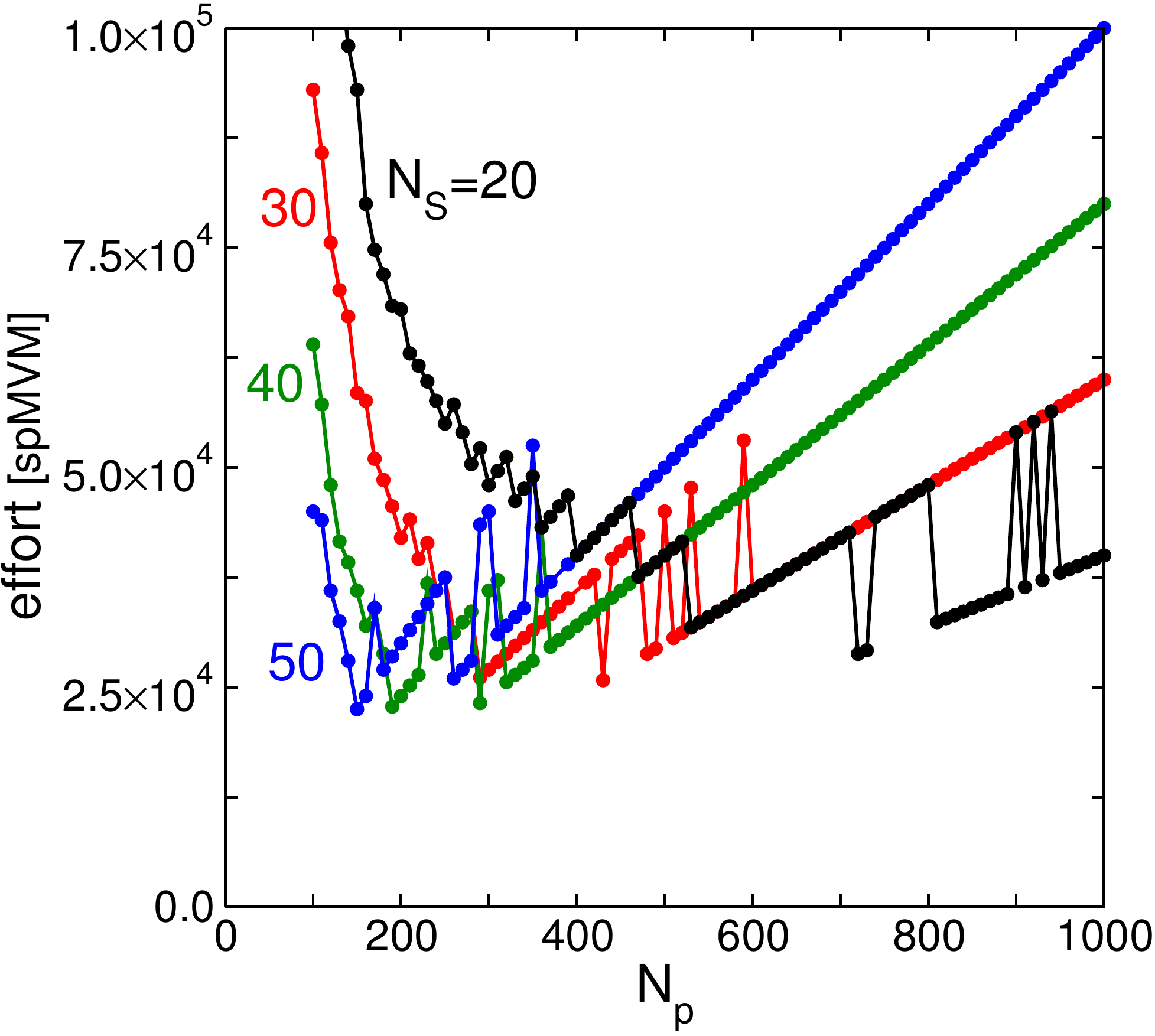}
\hspace*{\fill}
\caption{Effort of ChebFD, measured by the total number ($N_p \times N_S \times \textsf{iterations}$) of spMVMs until convergence, as a function of $N_S$ for fixed $N_p=100, 200, 500, 1000$ (left panel) and as a function of $N_p$ for fixed $N_S = 20, 30, 40, 50$ (right panel).
The matrix is the same diagonal matrix as used in Fig.~\ref{fig:filtermatrix},
with $D=10^3$ and target interval $[-0.01, 0.01] \subset [-1,1]$.
Convergence was detected according to the criterion stated in the text, with residual $\epsilon \le 10^{-7}$.
In all cases, the $N_T=10$ target eigenvalues were found with an accuracy better than $\epsilon$.
}
\label{fig:NumExp}
\end{figure}

An initial numerical experiment on the correct choice of $N_S$ and $N_p$ is shown in Fig.~\ref{fig:NumExp} for the diagonal matrix from Sec.~\ref{sec:Example1}.
The experiment clearly supports our previous advice:
Best convergence is achieved for $4 \le N_S/N_T \le 5$.
Smaller $N_S$ requires larger $N_p$, because the width of the filter polynomial must roughly match the width of the search interval.
In the present example with equidistant eigenvalues, i.e., a flat DOS, reducing $N_S \mapsto N_S/2$ should require increasing $N_p \mapsto 2 N_p$. This is consistent with the data in Fig.~\ref{fig:NumExp},
 until convergence becomes very slow for $N_S \lesssim 2 N_T$ for the reasons explained in Sec.~\ref{sec:Optimal}.

The present experiment can be summarized in a few rules of thumb:
First, the optimal choice of $N_p$ depends on $N_S$, and smaller $N_S$ requires a better filter polynomial with larger $N_p$.
Preferably, the number of search vectors is substantially larger than the number of target vectors, which increases the benefits of ``overpopulating'' the target space.
At least we should guarantee $N_S \ge 2 N_T$.
Second, the width of the central lobe of the filter polynomial,
which scales approximately as $N_p^{-1}$, should be of the same size as the search interval, hence, about twice as large as the target interval. Near the center of the spectrum $N_p \approx 4 (b-a)/(\lamo-\lamu)$ according to the plots\footnote{For comparison, see also footnote~\ref{afootnote}.} in Fig.~\ref{fig:filters}

These informal rules indicate what typically happens in a ChebFD run.
For less informal statements we have to return to the convergence theory from Secs.~\ref{sec:Theory} and \ref{sec:Optimal},
which requires information about the distribution of eigenvalues outside of the target interval, i.e., the DOS $\rho(\lambda)$.

\begin{table}[ht]
\centering
\small
\begin{tabular}{ccrrccc}
& & & & \multicolumn{2}{c}{effort $N_\mathrm{MVM}$ [$\times 10^6$]} \\
$N_S$ & $\delta'/\delta$ &  $N_p^\mathrm{opt}$ & \multicolumn{1}{c}{$\eta^\mathrm{opt}$} & estimate & numerical & iterations $N_\mathrm{iter}$ \\[0.5ex]\hline\hline
\multicolumn{7}{c}{flat DOS}\rule{0pt}{3.2ex} \\[0.5ex]
125 & 0.25 & 9972 & 4107 & 6.16 & 6.23 & 5\\
150 & 0.50 & 4988 & 2052 & 3.69 & 3.74 & 5 \\
175 & 0.75 &  3320 & 1368 & 2.87 &  2.90 & 5 \\
200 & 1.00 & 2500 & 1023 & 2.45 & 2.50 & 5\\
300 & 2.00 & 1258 & 507  & 1.83 & 2.30 & 6\rlap{$^*$}\\
400 & 3.00 &  817 & 351   & 1.68 & 1.96  & 6\\
500 & 4.00 &  612 & 282   & 1.69 & 2.14 & 7\rlap{$^*$}  \\
600 & 5.00 &  495 & 235   & 1.69 & 2.08 & 7\rlap{$^*$}  \\
\multicolumn{7}{c}{linear DOS}\rule{0pt}{3.2ex} \\[0.5ex]
125 & 0.12 & 1037 & 427  & 0.64 & 0.65 & 5 \\
150 & 0.22 & 565 & 233 & 0.42 & 0.42 & 5\\
175 & 0.32 & 388 & 160 &  0.34 &  0.34 & 5\\
200 &  0.41 & 303 & 125 & 0.30 & 0.30 & 5 \\
300 &  0.73 & 170 & 70 & 0.25 & 0.31 & 6   \\
400 &  1.00 & 124 & 51 & 0.24 &  0.30 & 6\\
500 &  1.24 & 99 & 41 & 0.25 & 0.30 & 6  \\
600 &  1.45 & 85 & 36 & 0.26 & 0.36 & 7\rlap{$^*$} \\
\end{tabular}
\caption{Numerical effort, in number of total spMVMs ($N_\mathrm{MVM}$), for the example of a flat and linear DOS. In both cases, the $N_T = 100$ central eigenvalues out of $D=40000$ eigenvalues are computed, with target interval $I_T=[-0.0025,0.0025]$ (flat DOS) or $I_T=[-0.05,0.05]$ (linear DOS),
and different number of search states $N_S$.
The accuracy goal is $\epsilon = 10^{-12}$.
Also reported are the search interval margin $\delta'$, the optimal polynomial degree $N_p^\mathrm{opt}$ and filter qualitity $\eta^\mathrm{opt}$, the theoretical estimate for $N_\mathrm{MVM}$ according to Eqs.~\eqref{etanflat},~\eqref{etanlinear}, and the number of iterations $N_\mathrm{iter}$ executed in the ChebFD scheme.
}
\label{tab:num}
\end{table}

\subsection{Numerical experiments for flat and linear DOS}

The considerations from Sec.~\ref{sec:Optimal} explain how the number of search vectors $N_S$ and the filter polynomial degree $N_p$ should be chosen depending on the DOS.
In Table~\ref{tab:num} we summarize numerical experiments for 
the computation of the $N_T = 100$ central eigenvalues of a $D=40000$ dimensional matrix.
For the flat DOS the target interval half width is $\delta =2.5 \times 10^{-3}$,
for the linear DOS $\delta = 0.05$.
The accuracy goal is $\epsilon = 10^{-12}$, the number of search vectors ranges from $
1.25 \le N_S/N_T \le 6$.
The estimate for the effort, measured by the total number of spMVMs, is $N_\mathrm{MVM} = \eta \times N_S \times (-\log_{10} \epsilon)$ as in Sec.~\ref{sec:Optimal}.
The numerical value is $N_\mathrm{MVM} = N_\mathrm{iter} N_S N_p$.

In agreement with the theoretical estimates from Eqs.~\eqref{etanflat},~\eqref{etanlinear} the numerical data show that the effort increases significantly for $N_S < 2 N_T$, becomes smaller for $N_S \gtrsim 2 N_T$, but does not decrease further for $N_S \gtrsim 4 N_T$.
Notice how the overall effort for the linear DOS is approximately smaller by the factor $4 \sqrt{N_T/D} = 1/5$ discussed previously.

In some cases (marked with a star in the table) the ChebFD algorithm executes one additional iteration in comparison to the theoretical estimate for $N_\mathrm{iter}$.
In these cases, the residuals in the previous iteration have already dropped to $\approx 2 \epsilon$ but not below $\epsilon$.
The one additional iteration increases the total effort by a large percentage 
if $N_S/N_T$ is too large. In practice, also the increasing time required for orthogonalization slows down the computation for very large $N_S$.
The data from Table~\ref{tab:num} thus confirm our recommendation to choose $N_S$ as a small multiple ($2 \le N_S/N_T \le 4$) of the number of target vectors $N_T$.

\section{Parallel implementation and performance engineering}
\label{sec:Impl}

The ChebFD scheme from the previous section is a straightforward and simple algorithm,
which allows for clean implementation and careful performance engineering.
The most work-intensive task, and hence the most relevant for performance analysis,  is step~\ref{corestep} shown in Fig.~\ref{fig:ChebFD}: the application of the polynomial filter.
A key feature of our implementation is the use of sparse matrix multiple-vector multiplication (spMMVM) as provided by the \ghost\ library~\cite{GHOST}, where the sparse matrix is applied simultaneously to several vectors.
As we demonstrated previously for the KPM~\cite{KHWPAF15} and a block Jacobi-Davidson algorithm~\cite{Roehrig14} the reduction of memory traffic in spMMVM can lead to significant performance gains over multiple independent spMVMs, where the matrix has to be reloaded from memory repeatedly.
Storage of block vectors in row major order is crucial to avoid scattered memory access and improve cache utilization.
Operations on block vectors are intrinsic to the ChebFD scheme (see the two
loops over $k$ in Fig.~\ref{fig:ChebFD}), and use of spMMVM promises substantial
performance gains. \ghost\ is an open-source library and is available for
download\footnote{\url{https://bitbucket.org/essex/ghost}}.

A thorough analysis of the optimization potential of block vector operations within the context of KPM is given in Ref.~\cite{KHWPAF15}.
The ChebFD computational core differs from KPM only in the addition of the \texttt{axpy} vector operations in lines 4,10 in Fig.~\ref{fig:ChebFD}.
The $N_S$ independent \texttt{axpy} operations in the outer loop can also be combined into a simultaneous block \texttt{axpy} operation with optimal cache utilization. 
In our implementation the \texttt{axpy} operations are incorporated into an augmented spMMVM kernel (contained in \ghost),
which saves $n_b$ vector (re)loads from memory during one combined spMMVM+\texttt{axpy} execution.
The intention behind the application of the optimization techniques \textit{kernel fusion} and
\textit{vector blocking} is to increase the computational intensity (amount of floating
point operations per transferred byte) of the original series of strongly
bandwidth-bound kernels including several BLAS-1 calls and the sparse matrix
vector multiply. 
As the number of floating point
operations is fixed for our algorithm, this can only be achieved by lowering the
amount of data transfers.

In our implementation we also compute the Chebyshev moments $\langle \vec x_k, \vec w_k \rangle$ during step~\ref{corestep},
which will allow us to monitor the spectral density contained in the current search space (as in Eq.~\eqref{weight1}) and thus check, e.g., the choice of $N_T$ during execution of the algorithm.
The computation of the dot products can be integrated into the spMMVM kernels
and does not affect  performance for kernels with low computational intensity because no additional memory accesses are
required.

Extending the analysis from Ref.~\cite{KHWPAF15}, we can adjust
the expression of the computational intensity (assuming perfect data re-use for
the spMMVM input block vector) $I(n_b)$ by taking the additional block
\texttt{axpy} operation into account. $S_d$ ($S_i$) denotes the size of a
single matrix/vector data (index) element and 
$F_a$ ($F_m$) indicates the number of floating point operations per addition
(multiplication).
$N_\mathrm{nzr}$ is the average
number of non-zero entries per matrix row. 
\begin{align}
    I(n_b) &= \frac{N_\mathrm{nzr}(F_a+F_m) + \lceil 9F_a/2 \rceil + \lceil 11F_m/2 \rceil}{N_\mathrm{nzr}/n_b(S_d+S_i) + 5S_d}\;\frac{\text{Flops}}{\text{Byte}} \\
    \intertext{Obviously, the parameters depend on the underlying system matrix
        and data types. Our computations are based on double precision data
        (complex for topological insulators and real for graphene) and
        $N_\mathrm{nzr}=13~(4)$ for topological insulators (graphene). Thus, for our application scenarios the following
    computational intensities of the compute kernel can be given:}
    I(n_b)^\text{Topi} &= \frac{146}{260/n_b+80}\;\frac{\text{Flops}}{\text{Byte}} \;, \\
    I(n_b)^\text{Graphene} &= \frac{19}{48/n_b+40}\;\frac{\text{Flops}}{\text{Byte}} \;.
\end{align}

\subsection{Hardware testbed}
We analyze the performance on an Intel Xeon E5-2697v3 (``Haswell-EP'')
CPU. It features 14 cores running at a nominal clock frequency of 2.6 GHz.
This CPU implements the AVX2 instruction set which involves 256-bit wide vector
processing units and fused multiply-add instructions (FMA). Hence, the
theoretical double precision floating point peak performance $P^\text{peak}$ adds up to 582.4
Gflop/s. The maximum attainable memory bandwidth $b$ as measured with a read-only
micro-benchmark (vector reduction) is $b=64.7$ GB/s. This CPU is installed in the
second phase of the SuperMUC\footnote{\url{https://www.lrz.de/services/compute/supermuc/}}
cluster at the Leibniz Supercomputing Centre (LRZ) on which the large-scale
experiments of this work have been performed. The Intel Compiler version 15.0
has been used throughout all experiments, and Simultaneous Multithreading
(SMT) threads have been disabled.

The considered CPU has the \emph{Cluster on Die} (CoD) operating mode enabled.
This novel feature splits up a multicore CPU socket
into two equally sized non-uniform memory access (NUMA) domains, which should
increase the performance for workloads which are NUMA-aware and bound to the
last level cache or main memory. 
However, the number of NUMA domains doubles within the compute nodes which may
lead to significant performance degredation for implementations which are not
NUMA-aware or bandwidth limited applications with a dynamic workload scheduling
strategy.

\subsection{Performance modelling and analysis}
The performance modeling approach largely follows the work presented in
Ref.~\cite{KHWPAF15}
where we focused on an Intel Xeon Ivy Bridge CPU. In that work we have demonstrated that the performance is
limited by main memory bandwidth for small $n_b$ and by on-chip execution
for intermediate to large $n_b$.
In the following we apply the same procedure to the Haswell-EP architecture,
which provides 29\% higher main memory bandwidth and a twice higher
floating point performance due to the aforementioned FMA instructions.

Assuming memory-bound execution, the roof{}line performance
model~\cite{Williams09} predicts an upper performance bound $P^*$
as the product of computational intensity and attainable memory bandwidth:
\begin{equation}
    P^* = I(n_b)b \;.
\end{equation}
As mentioned above, our formula for the computational intensity assumes perfect
data re-use of the spMMVM input block vector. Due to the limited cache size, this 
assumption is likely to be violated as the number of vectors $n_b$ increases.
The resulting data traffic overhead can be quantified as $\Omega =
V_\text{meas}/V_\text{min} \geq 1$ with $V_\text{min}$ ($V_\text{meas}$)
being the minimum (measured)
data volume. The data traffic measurements were done with hardware
performance counters using LIKWID~\cite{Treibig10}. 
This leads to a more accurate performance bound of
$P^*/\Omega$. 

\begin{figure}[tbp]
\begin{minipage}{0.45\textwidth}
\includegraphics[width=1\textwidth]{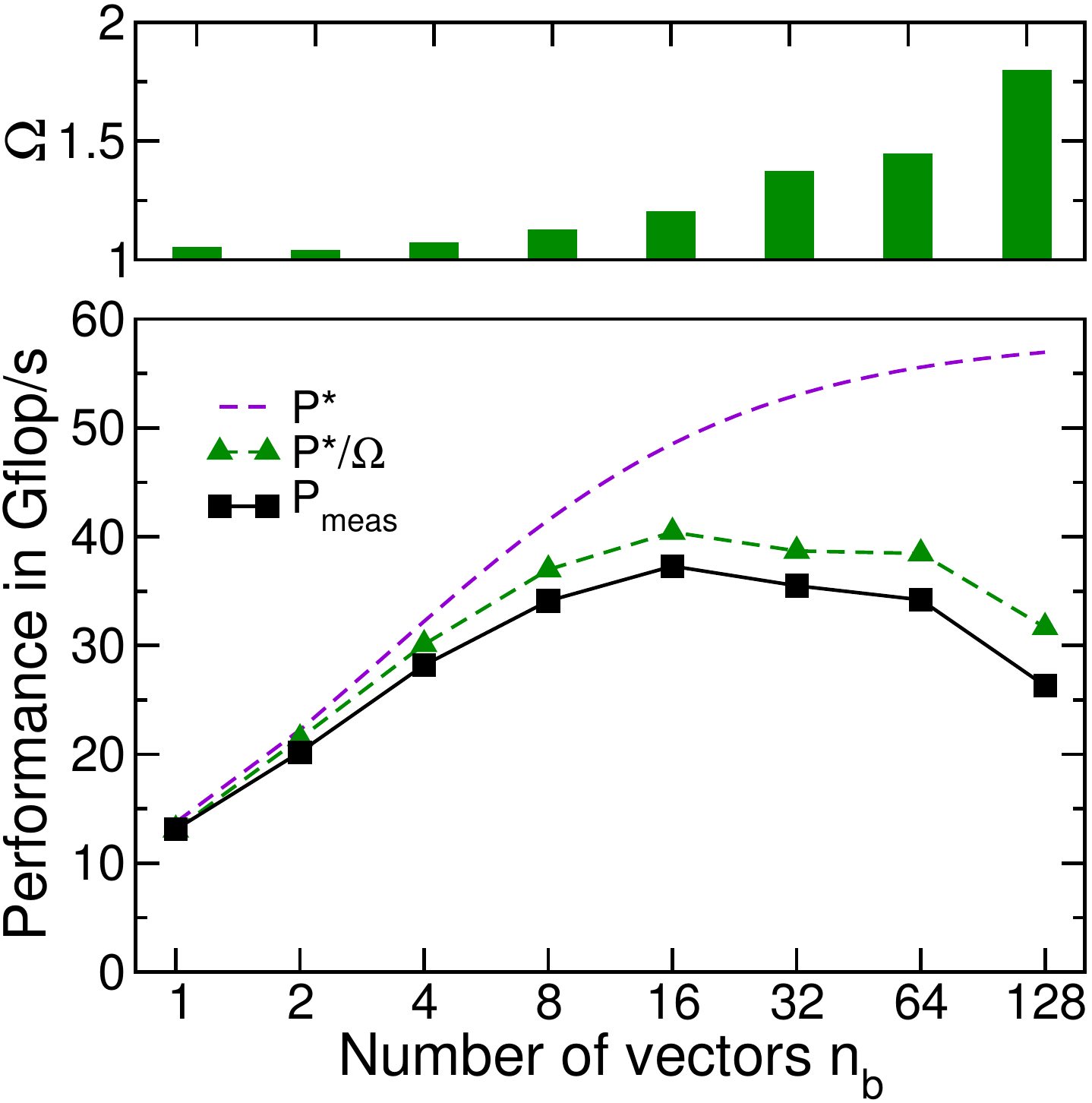}
\end{minipage}
\hspace*{\fill}
\begin{minipage}{0.45\textwidth}
\includegraphics[width=1\textwidth]{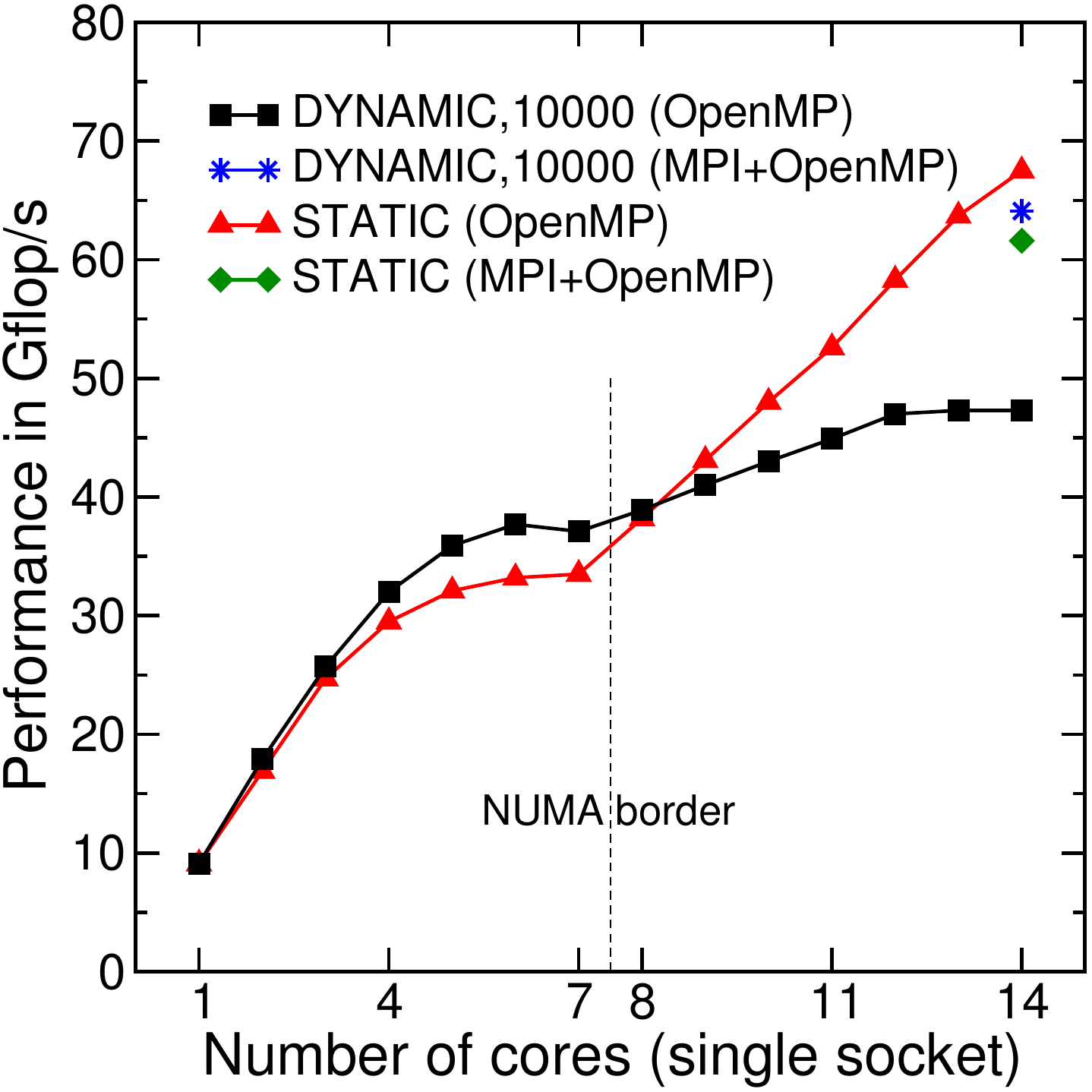}
\end{minipage}
\caption{
 Performance and data traffic overhead $\Omega$ (left panel) for the polynomial filter kernel
 (Fig.~\ref{fig:ChebFD}) as a function of the block size $n_b$ in a single
 NUMA domain (7 cores, pure OpenMP) with a topological insulator matrix for problem size $128 \times  64 \times 64$
 (matrix dimension $D =2^{21}$). In order to guarantee vectorized execution, the
 matrix was stored in SELL-32-1 (with \texttt{DYNAMIC,1000} OpenMP scheduling)
 for $n_b=1$ and SELL-1-1 (\texttt{DYNAMIC,10000}) otherwise. On the right
 panel we show scaling behavior inside a NUMA domain and in a
 full socket at $n_b=16$ for different OpenMP scheduling strategies
 and MPI+OpenMP hybrid execution (one process per NUMA domain).
 }
\label{fig:single_socket_perf}
\end{figure}

It can be seen in Fig.~\ref{fig:single_socket_perf} that the
measured performance $P_\text{meas}$ is always within 80\% of the prediction
$P^*/\Omega$, which substantiates the high efficiency of our pure OpenMP implementation in a
single NUMA domain.
Moreover we find that our implementation on Haswell is always limited by main
memory access and the on-chip execution no longer imposes a prominent bottleneck
(as it was the case for Ivy Bridge~\cite{KHWPAF15}). Thus, there is no need to implement further
low-level optimizations in our kernels. 
The deviation from the model increases with increasing $n_b$, which is caused
by other bottlenecks like in-core execution that gain relevance for the
increasingly computationally intensive kernel in this range. The best
performance in a single NUMA domain is obtained with \texttt{DYNAMIC} OpenMP
scheduling. However, this would lead to severe performance degradation on a full
socket due to non-local memory accesses. 
There are two ways to avoid this: Use \texttt{STATIC} OpenMP
scheduling or a hybrid MPI+OpenMP implementation with one MPI process per NUMA
domain. The performance of a statically scheduled OpenMP execution is shown 
in the right panel of Fig.~\ref{fig:single_socket_perf}. In a single NUMA domain
it is slightly lower
than for \texttt{DYNAMIC} scheduling, but thanks to the NUMA-aware
implementation in \ghost\ it scales across NUMA domains. Using
MPI+OpenMP instead of OpenMP only comes at a certain cost. For instance, input vector
data needs to be transferred manually which entails the assembly of
communication buffers. It turns out that the MPI+OpenMP performance with
\texttt{DYNAMIC} scheduling on a full socket is almost on par with the OpenMP-only
variant with \texttt{STATIC} scheduling; the MPI overhead is roughly compensated by the
higher single-NUMA domain performance. As expected, the statically scheduled
MPI+OpenMP variant shows the lowest performance on a full socket, as it unifies
both drawbacks: a lower single-NUMA domain performance and the MPI overhead. 
In summary for our application scenario there is no
significant
performance difference on a single socket between MPI+OpenMP+\texttt{DYNAMIC}
and OpenMP+\texttt{STATIC}.
In our large-scale experiments we chose the best variant on a single socket,
namely pure OpenMP with \texttt{STATIC} scheduling. 

\section{Large-scale application of ChebFD}
\label{sec:Large}

\subsection{Application scenario}

To demonstrate the potential of ChebFD for large-scale computations we choose an application scenario from current quantum physics research, the computation of the central eigenstates of a topological insulator as described in the introduction.
We use the model Hamiltonian~\cite{SRAF12}
\begin{equation}\label{Ham}
          H = -t \sum_{n,j=1,2,3} \left(  \Psi_{n+\hat \ee_j}^\dagger
                  \frac{\Gamma^1 - \ii \Gamma^{j+1} }{2}
                 \Psi_{n} + \text{H.c.}\right)  
               + \sum_{n} \Psi_{n}^\dagger \left( 
                      V_n     \Gamma^0
                    + 2       \Gamma^1
                   \right) {\Psi_{n}} \;,
\end{equation}
where the first term describes the hopping of electrons along the three-dimensional cubic lattice (from site $n$ to sites $n \pm \hat \ee_j$), and the second term contains the disorder potential $V_n$.
 Here, the $\Gamma^j$ are $4 \times 4$ Dirac matrices acting on the local $2 \times 2$ sublattice and spin degrees of freedom.
For a topological insulator slab of $L_x \times L_y \times L_z$ sites the matrix dimension is $D = 4 L_x L_y L_z$.

\begin{figure}
\hspace*{\fill}
\includegraphics[width=0.4\textwidth]{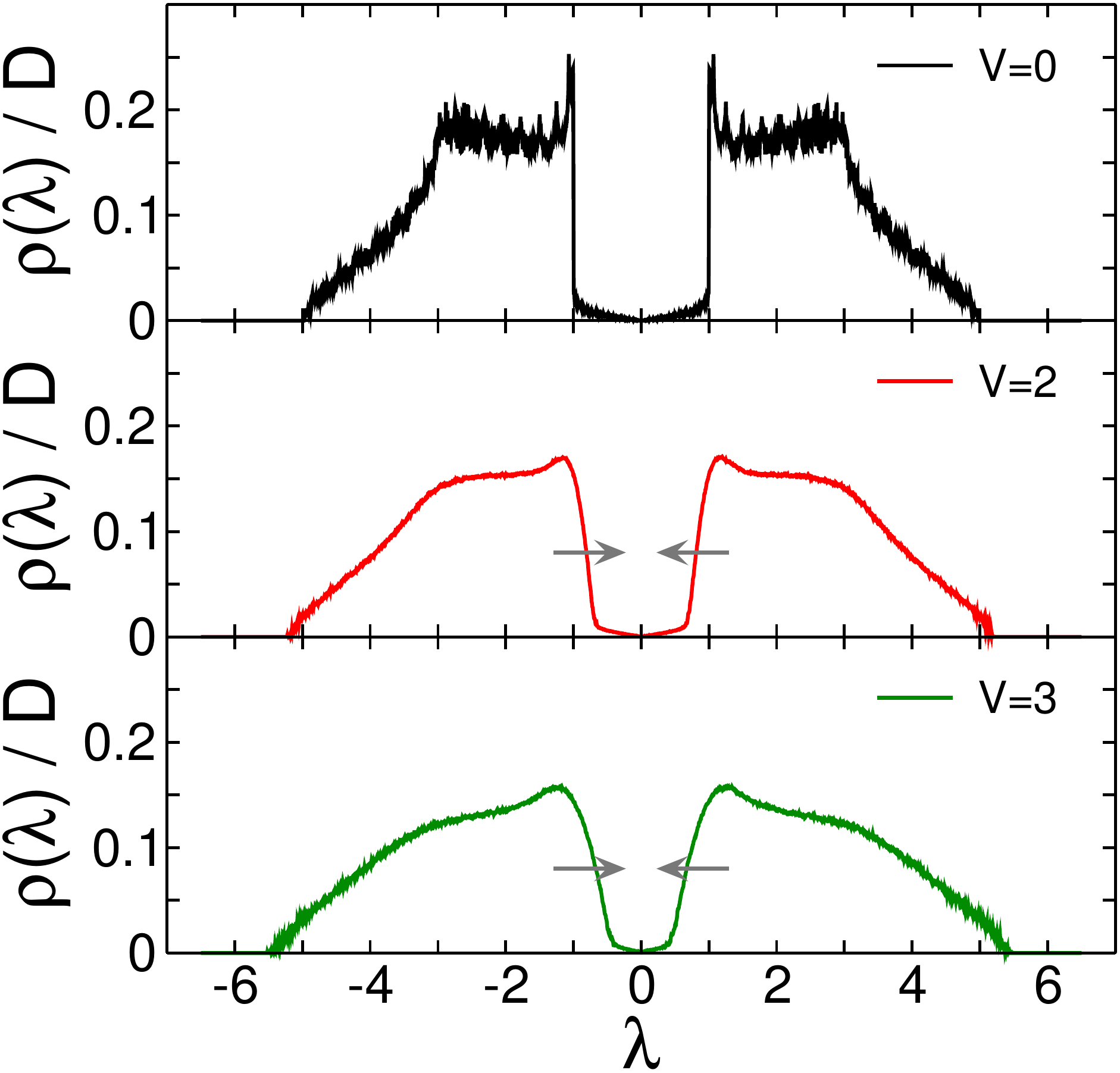}
\hspace*{\fill}
\includegraphics[width=0.4\textwidth]{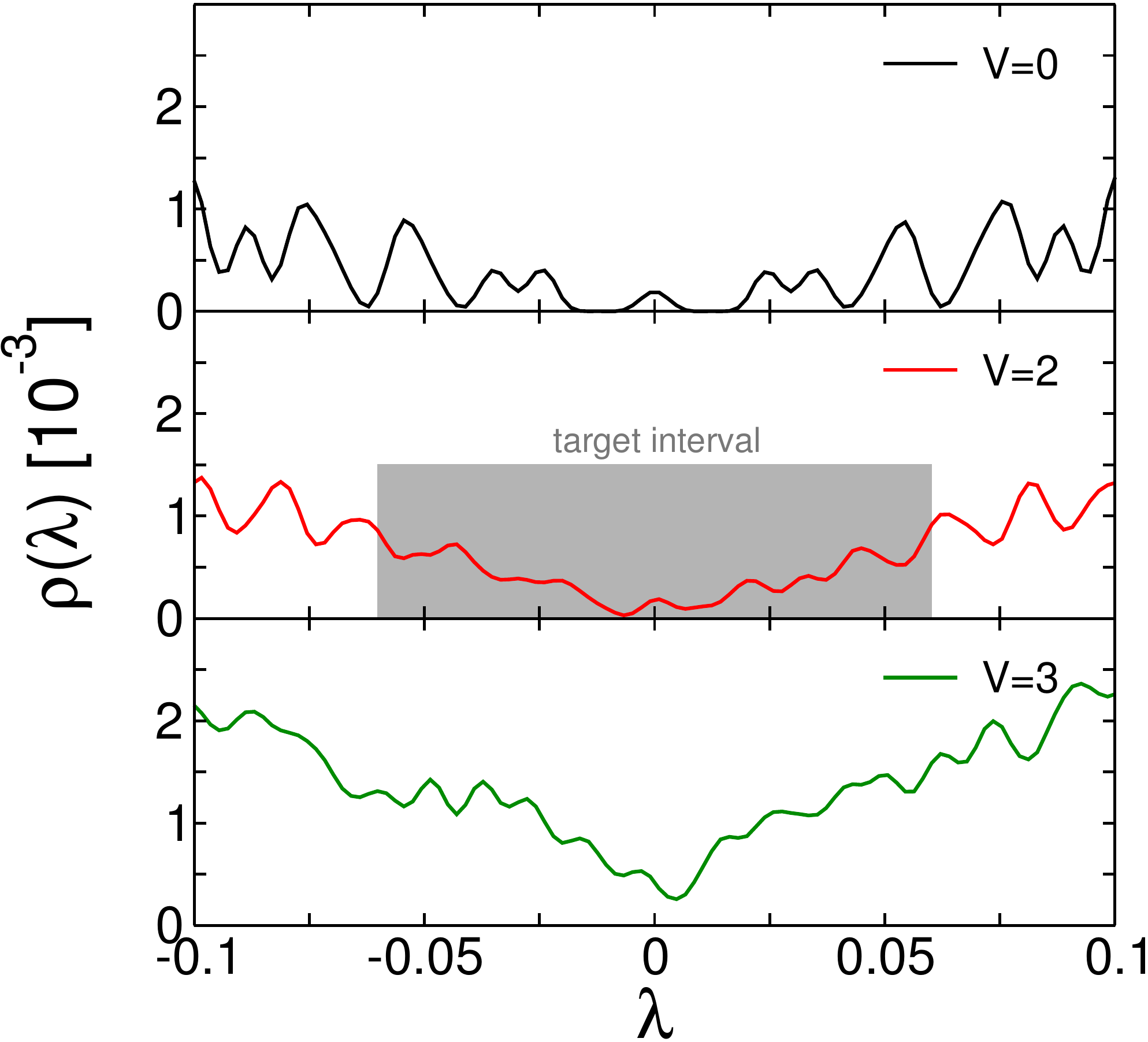}
\hspace*{\fill}
\caption{DOS $\rho(\lambda)$ for the topological insulator without ($V=0$) and with ($V \ne 0$) disorder, computed with KPM and stochastic trace evaluation~\cite{WWAF06}. 
Increasing the disorder shifts states into the gap at the center of the spectrum,
leading to the redistribution of the DOS as indicated by the arrows.
A zoom-in of the central region is shown in the right panels.
The target interval used for the ChebFD computations in the following figures is marked with the gray rectangle.
}
\label{fig:DOS}
\end{figure}

In Fig.~\ref{fig:DOS} we show the DOS of this model as a function of energy (i.e., eigenvalue $\lambda$ of $H$), without ($V_n \equiv 0$) and with ($V_n \in [-V/2, V/2]$, uniform box distribution) disorder.
Disorder fills the gap that exists in the DOS of the system with a small number of states.
Our goal is to compute these states, which are those most relevant for physics applications, with ChebFD.

\subsection{Numerical experiments}

\begin{figure}
\hspace*{\fill}
\includegraphics[width=0.45\textwidth]{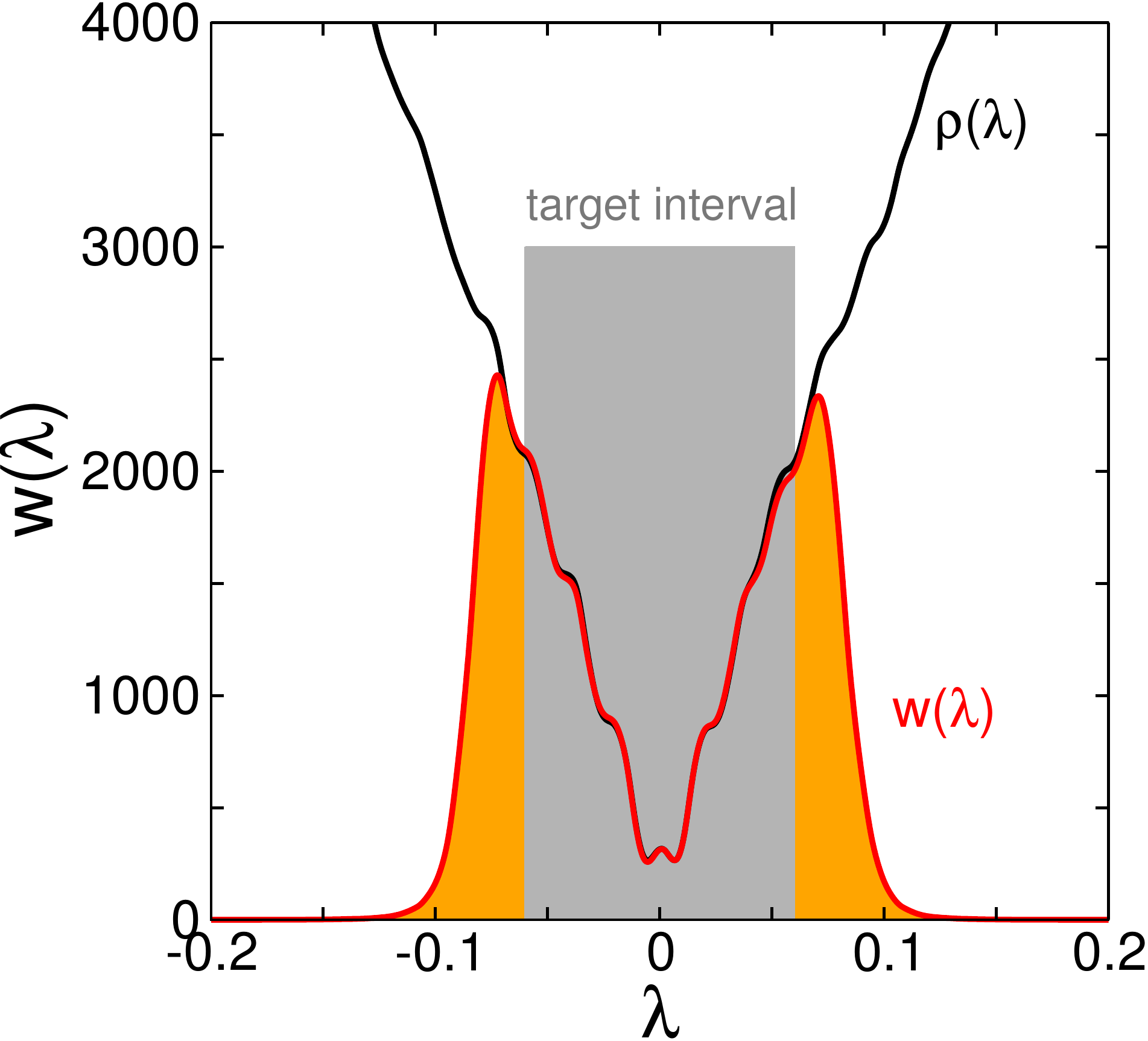}
\hspace*{\fill}
\includegraphics[width=0.45\textwidth]{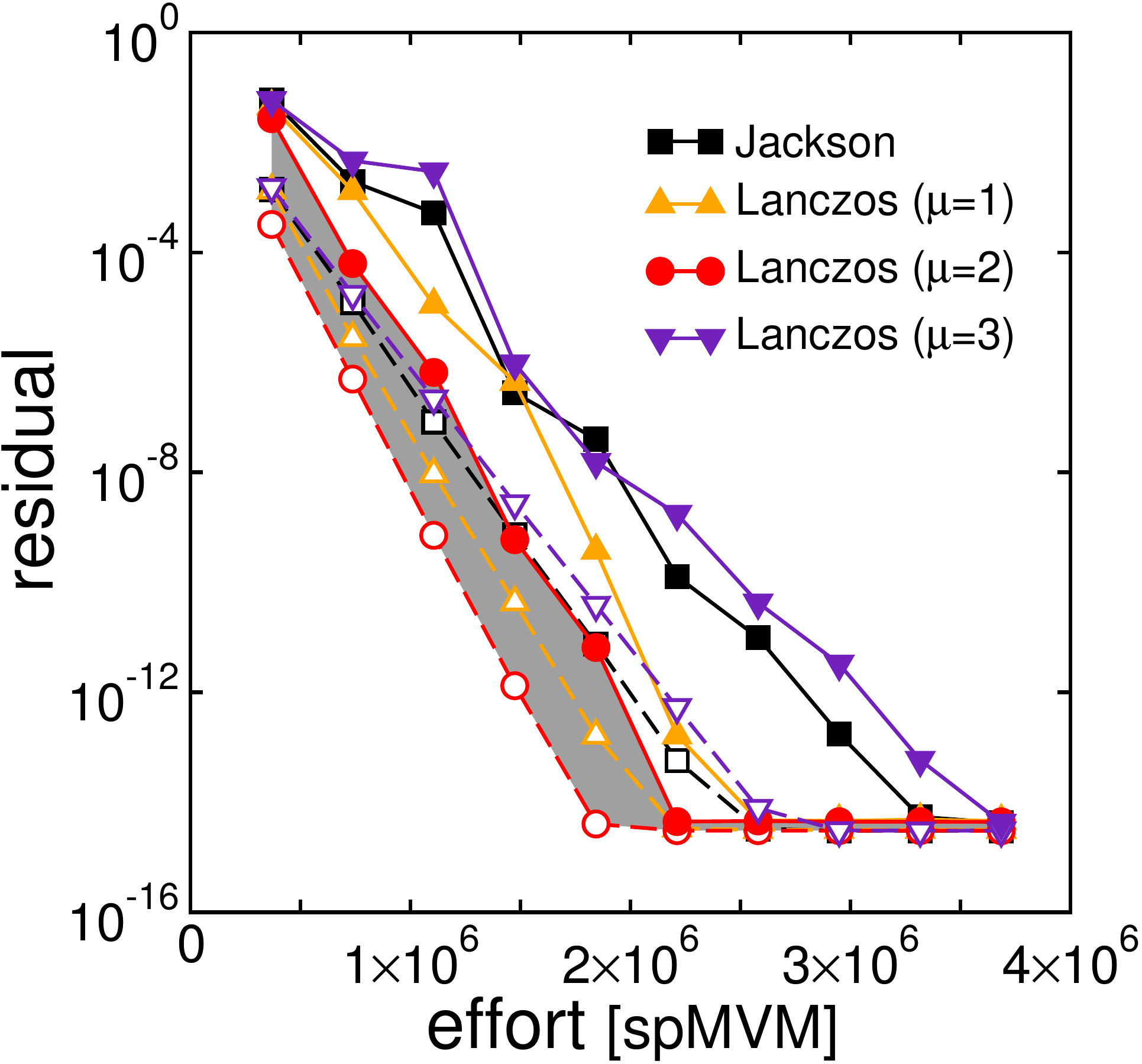}
\hspace*{\fill}
\caption{ Convergence of ChebFD with increasing number of iterations,
for the topological insulator example (Eq.~\eqref{Ham}, matrix dimension $2.6 \times 10^6$) with parameters as given in the text, and $N_S=256$, $N_p=1438$.
Left panel:
Weight $w(\lambda)$ of search vectors at target eigenvectors (see text) after 2 iterations. 
Shown is the DOS $\rho(\lambda)$ from Fig.~\ref{fig:DOS} for comparison.
Right panel:
Smallest (open symbols) and largest (filled symbols) residual
of the currently accepted Ritz vectors with Ritz values in the target interval.
Data are shown as a function of the total number ($N_p \times N_S \times \textsf{iterations}$) of spMVMs for the Jackson and Lanczos ($\mu = 1,2,3$) kernels.
The performance of the best kernel (Lanczos $\mu = 2$) is highlighted by the gray area.
}
\label{fig:ConvTopi}
\end{figure}

\begin{figure}
\hspace*{\fill}
\includegraphics[height=0.4\textwidth]{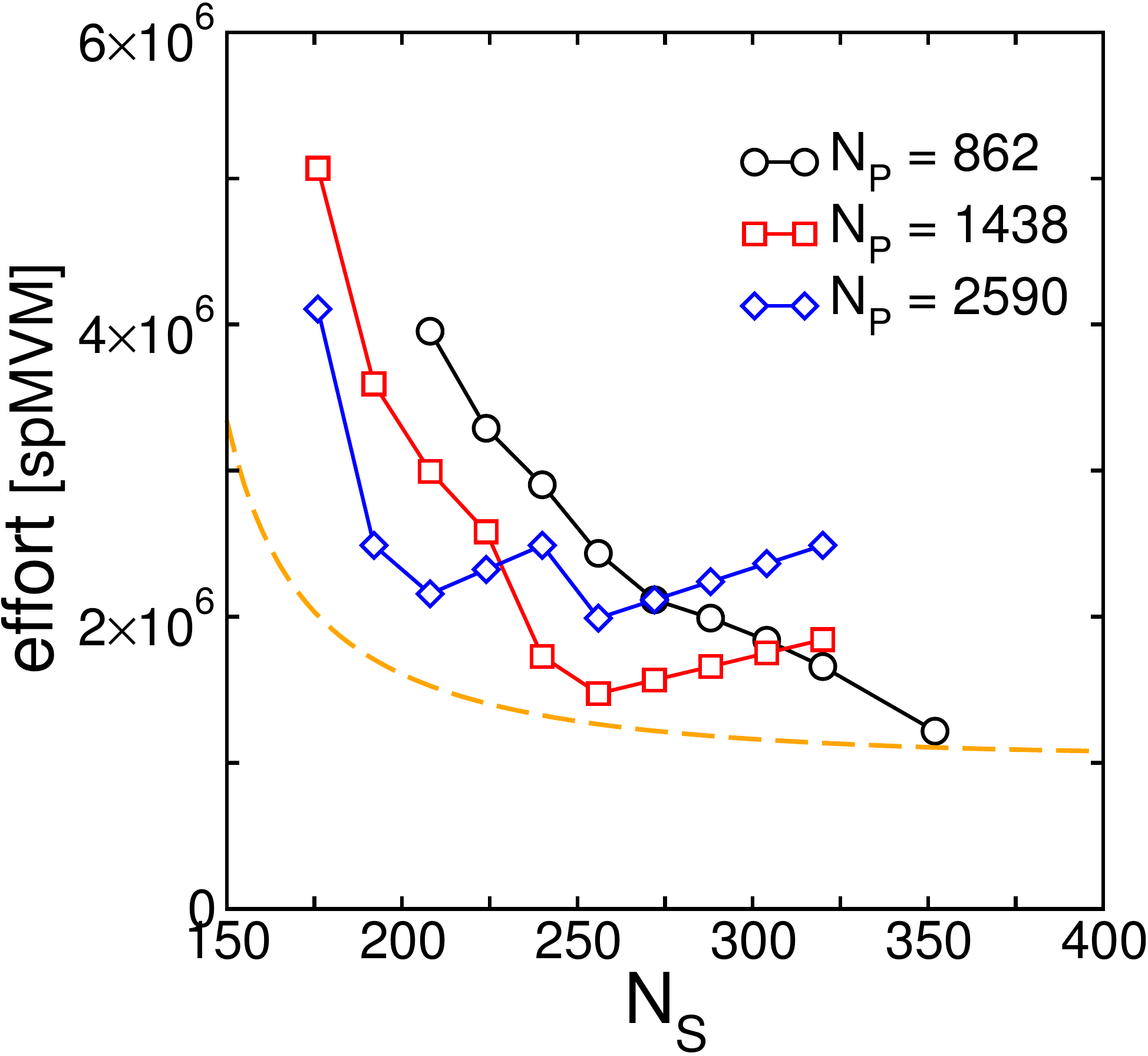}
\hspace*{\fill}
\includegraphics[height=0.4\textwidth]{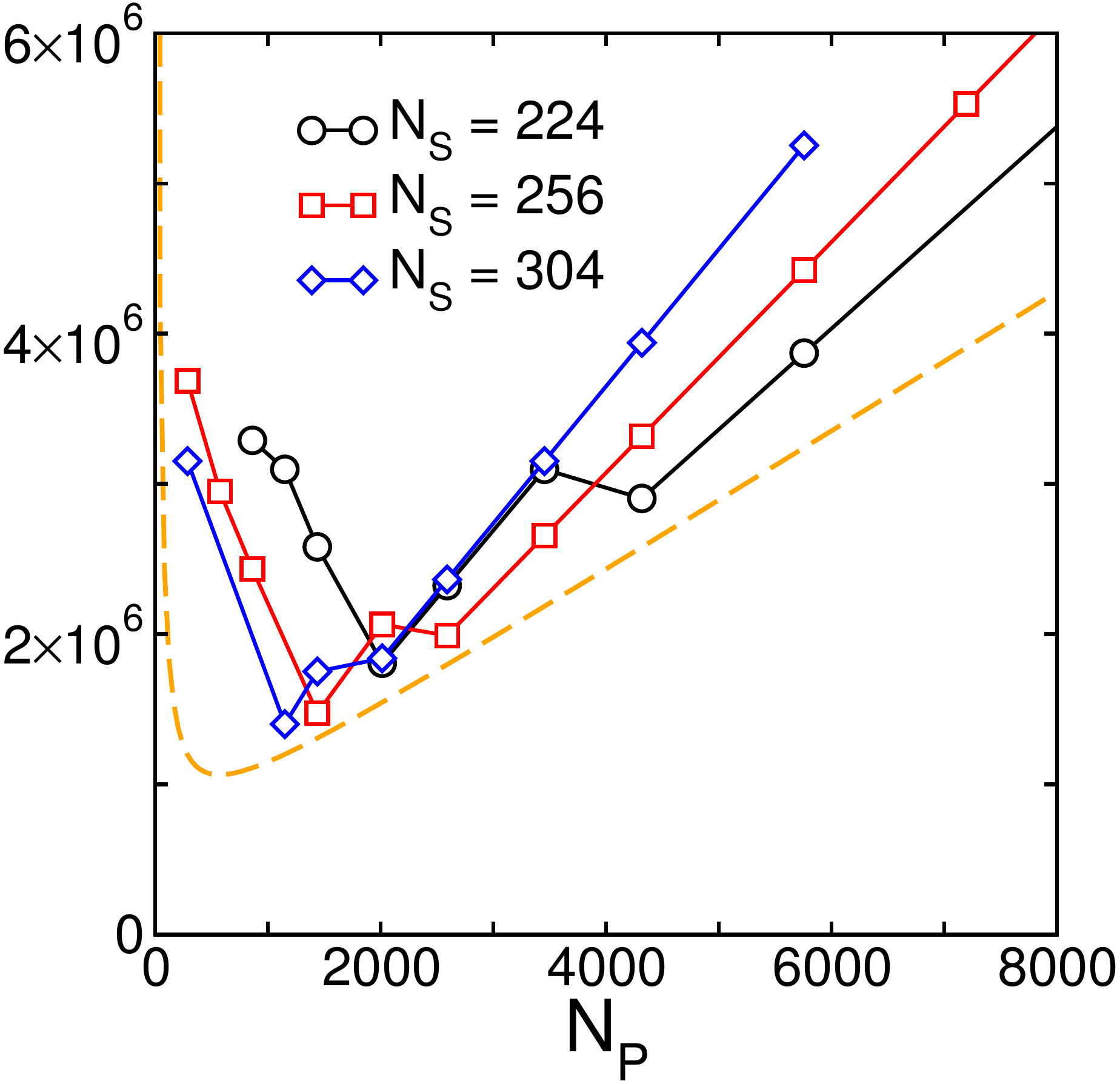}
\hspace*{\fill}
\caption{
Computational effort of ChebFD for different $N_S$ and $N_p$, similar to Fig.~\ref{fig:NumExp}, for the topological insulator example (Eq.~\eqref{Ham}) with parameters as given in the text.
The target space contains $N_T = 144$ vectors. The computation was stopped when all target vectors were found with residual below $\epsilon = 10^{-9}$.
Left panel: Number of spMVMs for fixed $N_p = 862, 1438, 2590$ as a function of $N_S$.
Right panel: Number of spMVMs for fixed $N_S=224, 256, 304$ as a function of $N_p$.
The dashed curves in both panels give the theoretical estimate of the minimal effort for a linear DOS (cf. Eq.~\eqref{etanlinear}).
\label{fig:EffortTopi}}
\end{figure}

Two initial numerical experiments are presented in Figs.~\ref{fig:ConvTopi} and \ref{fig:EffortTopi} for a topological insulator with size $256 \times 256 \times 10$ and matrix dimension $D = 5 \times 2^{19} \approx 2.6 \times 10^6$.
We set $t=1$, and the disorder strength is $V = 2.0$.
The matrix spectrum is contained in the interval $[a,b] = [-5.5, 5.5]$, the target interval is $[\lamu,\lamo]=[-0.06,0.06]$.
From the DOS $\rho(\lambda)$ in Fig.~\ref{fig:DOS} we estimate the number of target vectors as $N_T \approx 124$,
while the true value found after convergence of ChebFD is $N_T = 144$.
The difference (which is below $10^{-5} D$) is due to the stochastic sampling of the matrix trace.

In Fig.~\ref{fig:ConvTopi} we illustrate the convergence of ChebFD for this example.
In the left panel  we show, similar to Fig.~\ref{fig:filtermatrix}, the density $w(\lambda)$ of the (squared) weight $w_i$ of the search vectors at the target eigenvectors (cf. Eq.~\eqref{weight1}). 
Formally, $w(\lambda) = \sum_{i=1}^D w_i^2 \, \delta(\lambda - \lambda_i)$
with $\delta$-peaks at the eigenvalues $\lambda_i$ of $H$,
and we have $N_S = \int_a^b w(\lambda) d\lambda$.
Similar to the DOS we consider $w(\lambda)$ to be a continuous function for $N_S \gg 1$, rather than a distribution.
The quantity $w(\lambda)$ is conveniently estimated with KPM using the Chebyshev moments $\langle \vec x_k, T_H(H) \vec x_k \rangle = \langle \vec x_k, \vec w_k \rangle$, which can be computed during execution of the ChebFD core (step~\ref{corestep} from Fig.~\ref{fig:ChebFD}) without incurring a performance penalty.

We observe in Fig.~\ref{fig:ConvTopi} the same effect as discussed previously for Fig.~\ref{fig:filtermatrix}: search vectors are compressed into the target space through polynomial filtering.
Because the target space can accommodate only $N_T < N_S$ search vectors, finite weight accumulates also around the target interval.
As the residual of the Ritz pairs computed from the current search space decreases (right panel)
the search vector weight $w(\lambda)$ converges to the DOS $\rho(\lambda)$ inside of the target interval.
Note that also in this example `ghost' Ritz pairs with large residual occur in the target interval. Ritz pairs are accepted  as valid eigenpair approximations according to the criteria stated in Sec.~\ref{sec:ChebFD}.

A comparison of the different kernels (Jackson and Lanczos in Fig.~\ref{fig:ConvTopi}) shows that convergence is indeed fastest with the Lanczos kernel.
This result is consistently reproduced in other applications. Therefore, we 
recommend the Lanczos kernel and will use it in the remainder of the paper.

In Fig.~\ref{fig:EffortTopi} we repeat the numerical experiment from Fig.~\ref{fig:NumExp} for the present example. 
We observe again that (i) the computational effort increases rapidly for $N_S \lesssim 2 N_T$ (here, $N_S \lesssim 250$ in the left panel), and (ii) the computational effort becomes smaller if $N_S$ is increased to a small multiple of $N_T$.

The data in Fig.~\ref{fig:EffortTopi} fully agree with our analysis from Secs.~\ref{sec:Theory} and \ref{sec:Optimal}. 
Assuming a linear DOS in the vicinity of the target interval (which is justified by Fig.~\ref{fig:DOS})
the change of the computational effort with $N_S$ and $N_p$ can be predicted with the theoretical estimate from Eq.~\eqref{etanlinear} (dashed curves).

The computational effort is minimal for the largest number of search vectors used in the present example ($N_S = 350$ in the left panel). While the theoretical estimate predicts that the computational effort would be reduced further for larger $N_S \simeq 4 N_T \simeq 600$,  finite memory resources restrict the possible $N_S$ in large-scale applications.
The minimal effort predicted by theory (see text after Eq.~\eqref{etanlinear}) is $N_\mathrm{MVM} = 1.06 \times 10^6$ for a linear DOS, while our best computation required  $N_\mathrm{MVM} = 1.2 \times 10^6$. 
At the risk of overinterpreting this agreement (e.g., the DOS is not perfectly linear in the present example) we find that the algorithmic efficiency of our ChebFD implementation lies within $20 \%$ of the expected theoretical optimum.

\subsection{Large-scale topological insulator computations}

\begin{figure}
\hspace*{\fill}
\includegraphics[height=0.40\textwidth]{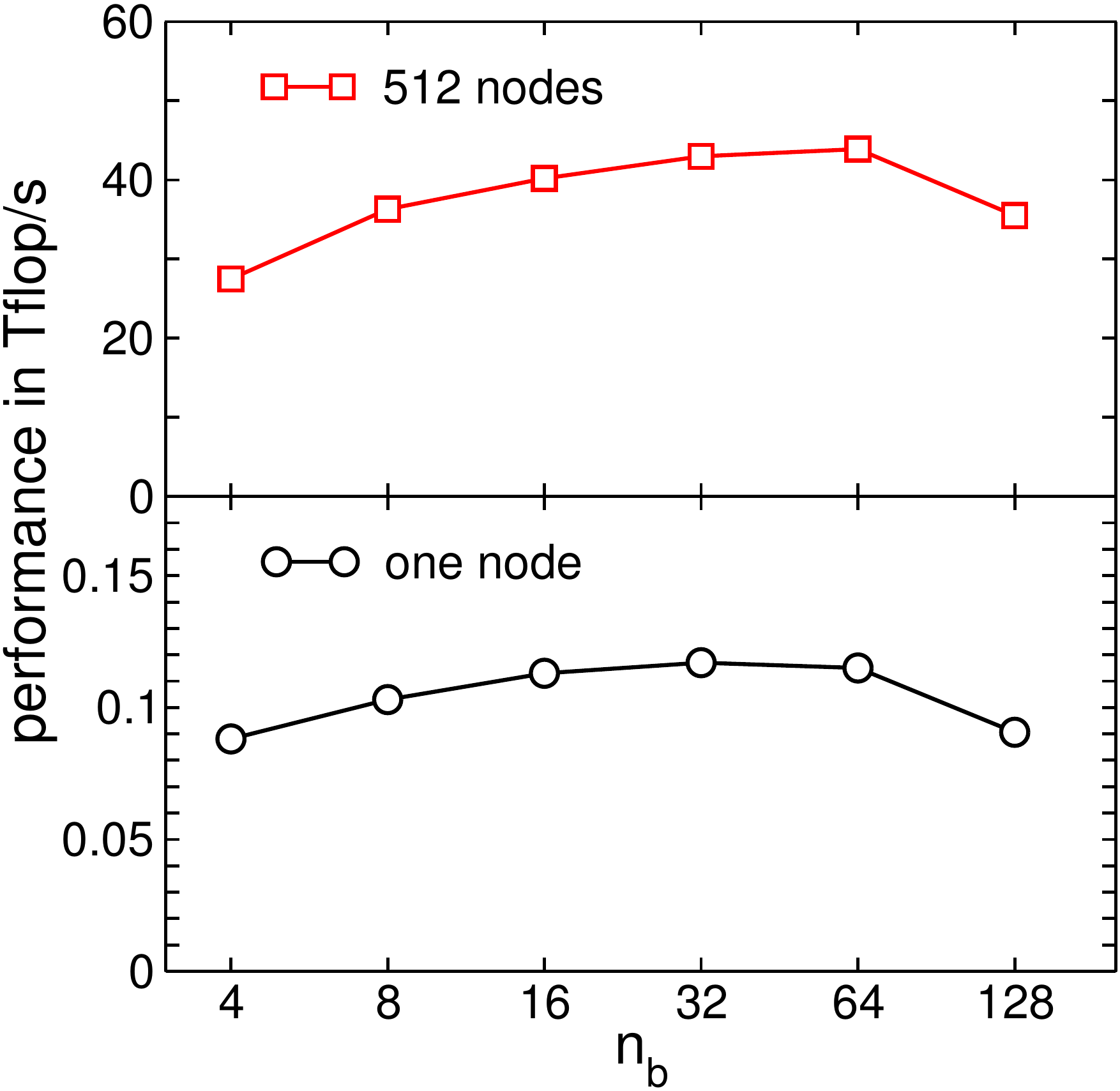}
\hspace*{\fill}
\includegraphics[height=0.40\textwidth]{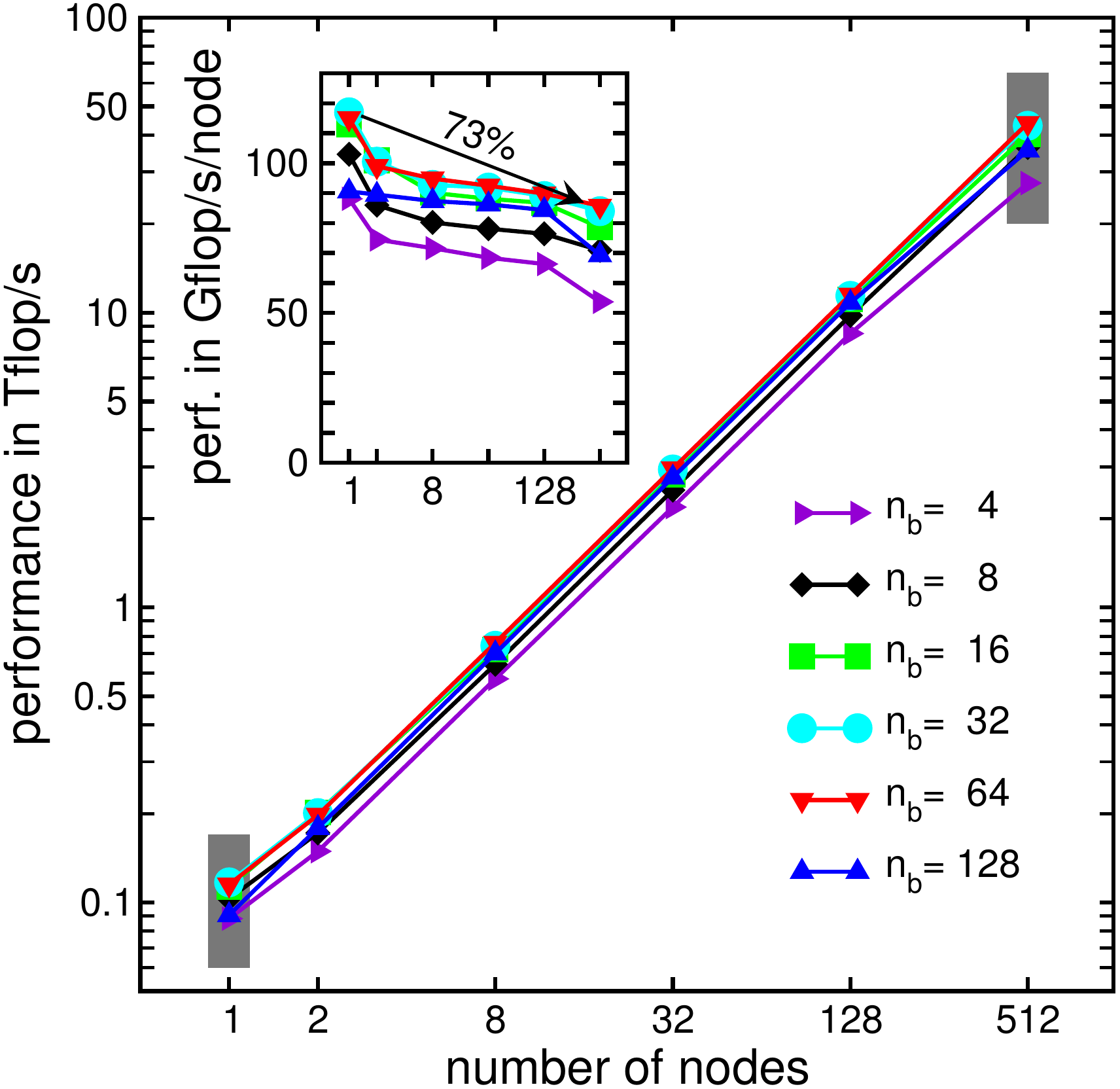}
\hspace*{\fill}
\caption{
Weak scaling performance for the polynomial filter kernel of the parallel ChebFD implementation,
for fixed problem size ($ 128 \times  64 \times 64$ slice, vector length $2^{21}$) per node.
Left panels: Sustained performance accumulated over all nodes as a function of the block size $n_b$ used in spMMVM, for a computation on 1 and 512 nodes.
Right panel: Sustained performance as a function of the number of compute nodes.
The data points highlighted in gray are those shown in the left panels.
The inset shows the average performance per node as an indicator for
the parallel efficiency.
\label{fig:Large}}
\end{figure}

We now examine the computational efficiency of our ChebFD implementation in large-scale computations.
For the topological insulator application,
where ChebFD is required for the solution of large problems that can not be treated with other techniques, the weak scaling performance is most relevant.

The weak scaling performance of the ChebFD computational core (step~\ref{corestep} from Fig.~\ref{fig:ChebFD}) is shown in Fig.~\ref{fig:Large},
together with the influence of the block size $n_b$ used in the spMMVM.
Computations took place on up to $512$ Haswell nodes of SuperMUC at LRZ Garching\footnote{see \url{https://www.lrz.de/services/compute/supermuc/systemdescription/} for the hardware description},
with fixed problem size ($D=2^{21} \approx 2 \times 10^6$) per node.
A correct choice of the block size $n_b$ can lead to a 2$\times$ gain in performance
(Fig.~\ref{fig:Large}, left panel),
which is in accordance with the single socket improvements demonstrated in Fig.~\ref{fig:single_socket_perf}.
Generally, moderately sized blocks (about $n_b = 2^6$) are preferential.
The overall scalability of ChebFD does not depend on $n_b$, as could be expected because spMMVMs optimize local memory traffic on the node-level but do not affect the inter-node communication volume.
However, larger block sizes become more effective with increasing number of nodes and problem size, because in our implementation the size of MPI messages during spMMVM is proportional to $n_b$.
Large $n_b$ thus reduces the inter-node communication effort in terms of latency and message count.
Comparing the single-node performance numbers in Fig.~\ref{fig:Large} to the single-socket performance as
discussed in Sec.~\ref{sec:Impl}, we can observe a significant drop of parallel
efficiency when going from one to two sockets in a single node, which is due to
the enabling of MPI communication.
Beyond a single node we achieve a parallel efficiency of 73\% when scaling from 1 to 512 nodes. 
In Fig.~\ref{fig:Large} two drops of parallel efficiency can be observed in the
scaling range:
from one to two nodes we are effectively enabling
communication over the network and from 128 to 512 nodes we encounter an effect
of the special network topology of the SuperMUC cluster.
It is built of islands of 512 nodes which internally have a fully
non-blocking FDR14 Infiniband network. However, a pruned tree connects
those islands, causing the intra-island bandwidth to be a factor of four larger
than inter-island. While all jobs up to 128 nodes ran in a single island,
512-node jobs were spread across two island by the batch system.
Both the step from one to two nodes and the step from one to two islands lead to an increased influence of
communication and to a decrease of parallel efficiency.

\subsection{Large-scale benchmark runs}

Additional benchmark data  are shown in Tab.~\ref{tab:benchmark},
using topological insulator and graphene matrices with dimensions ranging up to $D \gtrsim 10^9$.
The graphene matrices are obtained from the standard tight-binding Hamiltonian on the honeycomb lattice with an on-site disorder term~\cite{CGPNG09}. Generator routines for all matrix examples are provided with the \ghost\ library.

At present, interior eigenvalues cannot be computed for such large matrix dimensions  with methods based on (direct or iterative) linear solvers, and thus require a polynomial filtering technique like ChebFD.
One observes that with increasing matrix dimension the evaluation of the filter polynomial constitutes an increasing fraction of the total runtime and quickly becomes the dominating step.
Our optimizations for the ChebFD computational core are thus especially relevant for large-scale computations, and the benchmark data from Fig.~\ref{fig:Large} and Tab.~\ref{tab:benchmark} clearly demonstrate the potential of our ChebFD implementation.

Note the very high polynomial degree $N_p$ required for the graphene matrices,
where the target interval comprises less than $0.1\%$ of the entire spectrum.
Evaluating filter polynomials causes no difficulty even for such extremely high degrees because of the inherent stability of the Chebyshev recurrence.

\begin{table}
\begin{tabular}{c|ccccrrr}
  matrix & nodes & $D$ & $[\lamu : \lamo]_\mathrm{rel} $ & $N_T$   & $N_p$ & \parbox[b]{1.25cm}{runtime\\{}  [hours]} & \parbox[b]{1.75cm}{sust. perf.\\{} [Tflop/s]} \\ \hline
\texttt{topi1}     &  32 &  6.71e7 &   7.27e-3   & 148 &  2159 &  3.2  (83\%) & 2.96\hspace*{2em} \\
\texttt{topi2}     & 128 &  2.68e8 &   3.64e-3   & 148 &  4319 &  4.9  (88\%)  & 11.5\hspace*{2em} \\
\texttt{topi3}     & 512 &  1.07e9 &   1.82e-3   & 148 &  8639 & 10.1  (90\%) & 43.9\hspace*{2em} \\
\texttt{graphene1} & 128 &  2.68e8 & 4.84e-4 & 104 & 32463 & 10.8  (98\%) & 4.6\hspace*{2em} \\
\texttt{graphene2} & 512 &  1.07e9 & 2.42e-4 & 104 & 64926 & 16.4  (99\%) & 18.2\hspace*{2em} 
\end{tabular}
\caption{
ChebFD benchmark runs on up to 512 SuperMUC Haswell nodes for topological insulator (\texttt{topi\_}) and graphene (\texttt{graphene\_}) matrices,
 with matrix dimension $D$ up to $10^9$.
The central columns give the relative width of the target interval $[\lamu,\lamo]$ in comparison to the entire spectrum, the number of computed target eigenpairs ($N_T$) and the degree of the filter polynomial ($N_p$).
The last two columns give the total runtime
with the percentage spent in the polynomial filter evaluation step from Fig.~\ref{fig:ChebFD} and the sustained performance of this step.
We used $N_S = 256$ for all examples, and $n_b = 64$ ($128$) for the topological insulator (graphene) matrices.
}
\label{tab:benchmark}
\end{table}

\subsection{Topological insulator application example}

\begin{figure}
\hspace*{\fill}
\raisebox{-0.5\height}{\includegraphics[width=0.24\textwidth]{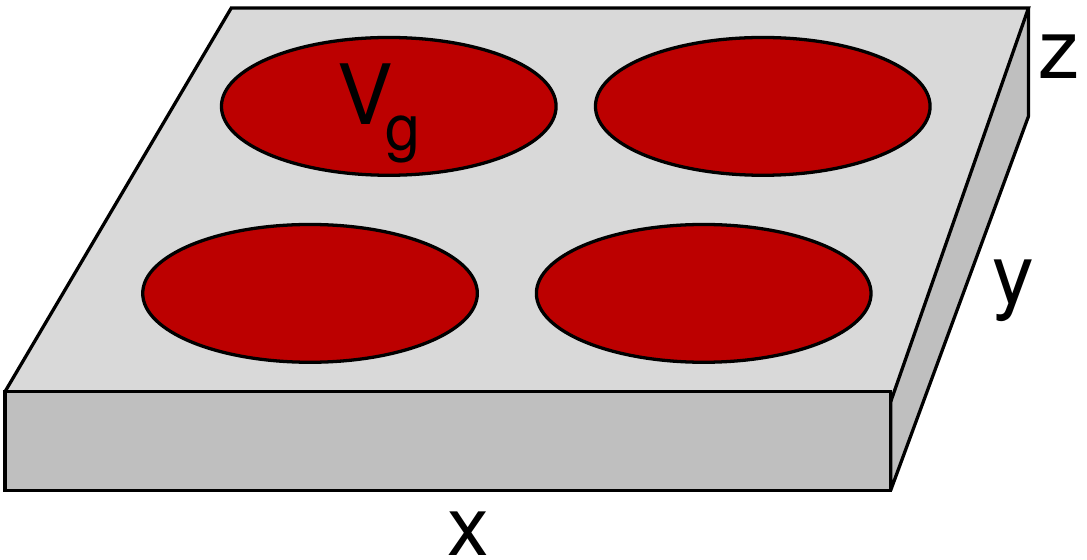}}
\hspace*{\fill}
\raisebox{-0.5\height}{\includegraphics[width=0.4\textwidth]{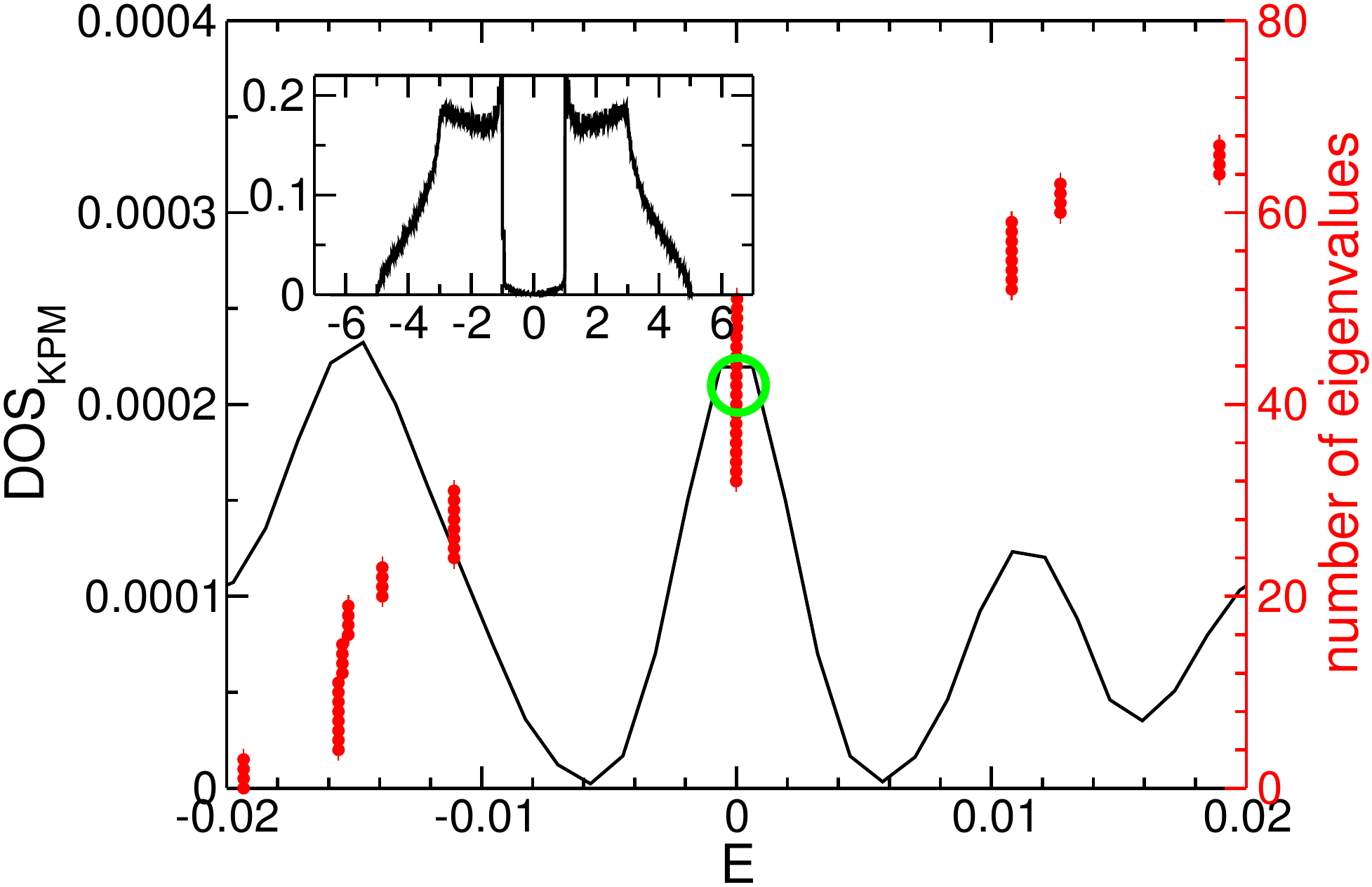}}
\hspace*{\fill}
\raisebox{-0.5\height}{\includegraphics[width=0.3\textwidth]{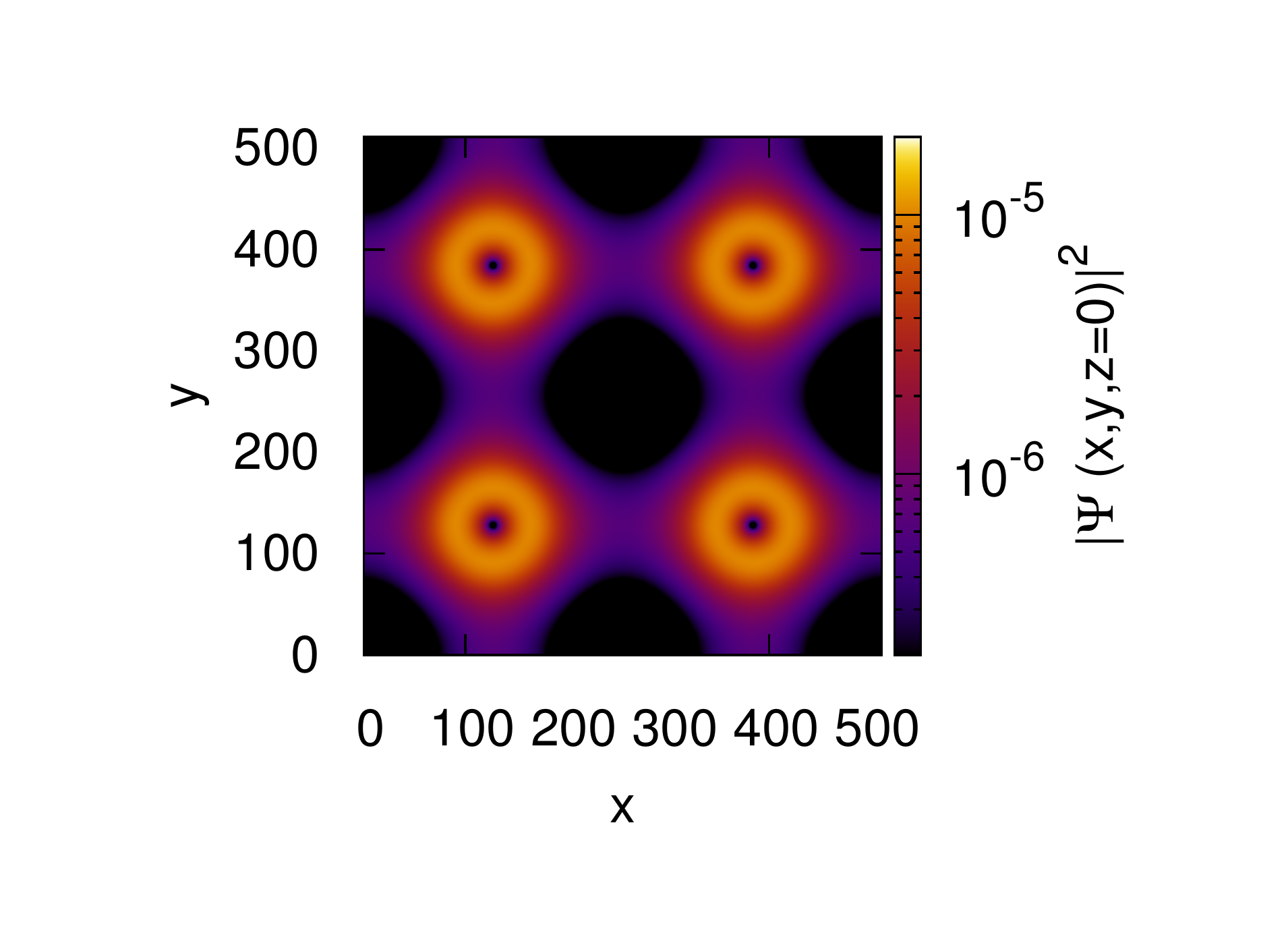}}
\hspace*{\fill}
\caption{Schematic geometry of gate-defined quantum dots in a topological insulator sheet (left panel), DOS $\rho(\lambda)$ and interior eigenvalues computed with ChebFD (central panel), and the wave function $|\psi(x,y,z=0)|^2$ in the $xy$ plane 
belonging to the eigenvalue at the Fermi energy marked with the green circle (right panel).}
\label{fig:Dot1}
\end{figure}

The good scalability and high performance of our ChebFD implementation allows us to deal with the very large matrices that arise in the theoretical modeling of realistic systems.
To illustrate the physics that becomes accessible through large-scale interior eigenvalue computations we consider the example of gate-defined quantum dots in a topological insulator (Fig.~\ref{fig:Dot1}), a primary example for the functionalization of topological materials.

In this example a finite gate potential is tuned to a value ($V_g= 0.11974$) such that the wave functions of the innermost eigenstates, i.e., those located at the Fermi energy, are concentrated in the circular gate regions. 
Minimal changes of the gate potential can lead to large changes of the wave functions, which gives rise to a transistor effect that allows for the rapid switching between an `on' and an `off' state.
While this effect can be studied already for small systems, or is even accessible to analytical treatment due to the high symmetry of the situation, the modeling of realistic systems requires the incorporation of disorder at the microscopic level, which results either from material imperfections or deliberate doping with impurity atoms.

\begin{figure}[hb]
\hspace*{\fill}
\raisebox{-0.5\height}{\includegraphics[width=0.48\textwidth]{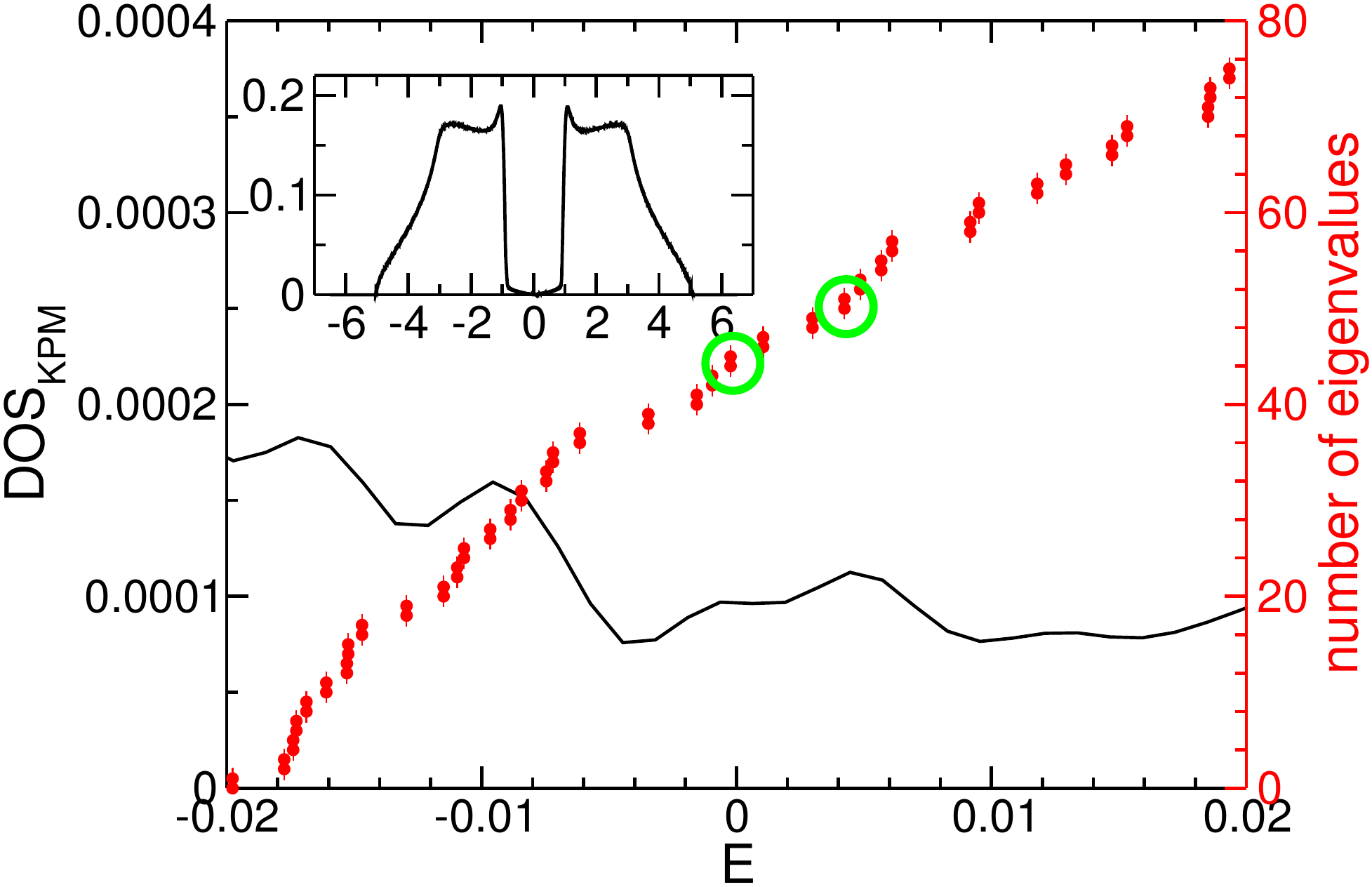}}
\hspace*{\fill}
\begin{minipage}[c]{0.38\textwidth}
\includegraphics[width=0.8\textwidth]{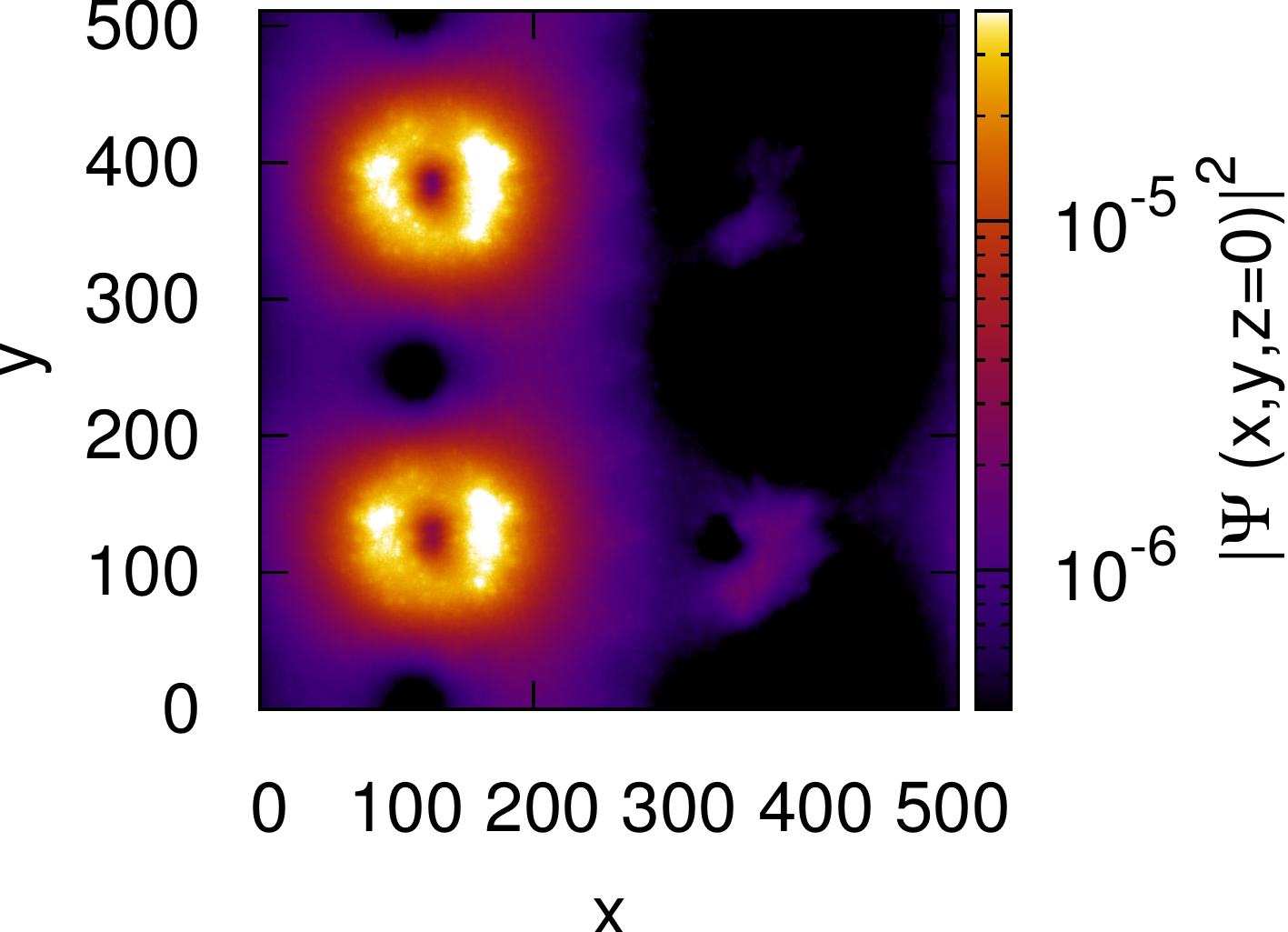} \\[+1ex]
\includegraphics[width=0.8\textwidth]{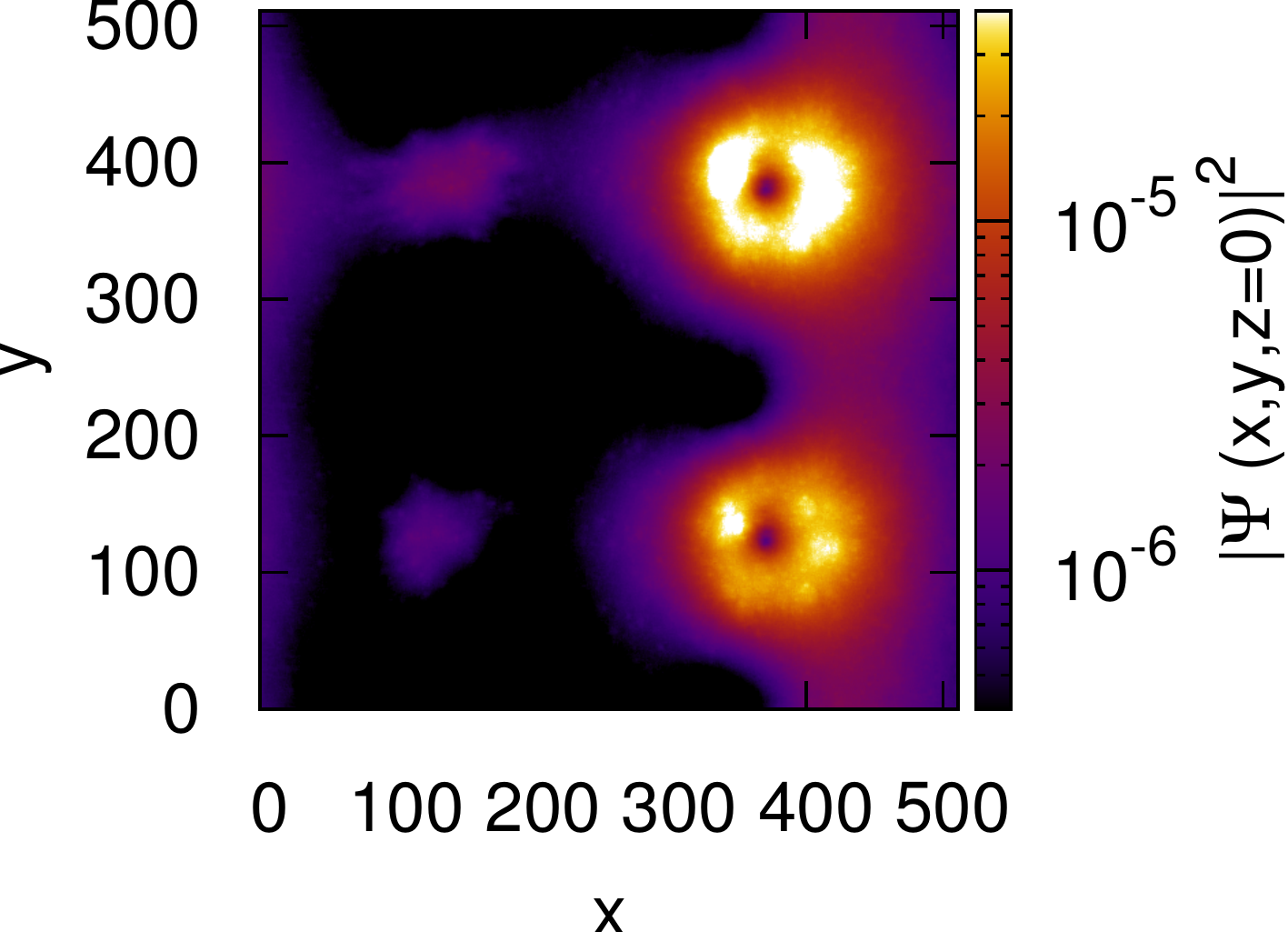}
\end{minipage}
\hspace*{\fill}
\caption{Left panel: DOS and interior target eigenvalues computed with ChebFD for a disordered (impurity-doped) topological insulator sheet with gate-defined quantum dots as in the previous figure.
The system size is $L_x \times L_y \times L_z = 2^8 \times 2^8 \times 2^4$, i.e., the matrix dimension is $D = 2^{26} \simeq 1.6 \times 10^7$, with disorder $V=1$.
Right panel: Wave function $|\psi(x,y,z=0)|^2$ in the $xy$ plane belonging to the two eigenvalues marked by green circles in the left panel.}
\label{fig:Dot2}
\end{figure}

The natural question then is how disorder modifies the transistor effect.
To answer this question the interior eigenstates must now be computed numerically, and large system sizes are required to correctly describe the effect of the microscopic disorder on the mesoscopic quantum dot wave functions.
An example for such a computation with ChebFD is shown in Fig.~\ref{fig:Dot2}.
Eigenstates close to the Fermi energy are still concentrated in the gate regions, although the circular shapes visible in the wave function from Fig.~\ref{fig:Dot1} are washed out by disorder.
Interestingly, the two nearby eigenstates shown in the picture are concentrated in different gate regions, which is possible because disorder breaks the perfect square symmetry of the arrangement.
Therefore, disorder can induce a new transistor effect,
where a small gate potential switches an electric current between different sides of the sample.
Full exploration of this effect certainly requires more realistic modeling of the situation.
As the present examples demonstrate, ChebFD can become a relevant computational tool for such explorations.

\section{Conclusions}
\label{sec:Conclude}

Polynomial filtering techniques such as the ChebFD scheme examined in the present paper are powerful tools for large-scale computations of many interior eigenvalues of sparse symmetric matrices.
These techniques combine two aspects.
On the one hand, they use polynomial filter functions to damp out eigenvector contributions outside of the target interval.
Polynomial filters allow for simple evaluation through spMVM, and are thus applicable to large matrices.
On the other hand, these techniques rely on large search spaces,
which can lead to fast convergence also for moderate quality of the polynomial filters.

ChebFD is a straightforward implementation of polynomial filtering.
Despite its algorithmic simplicity, ChebFD is well suited for large-scale computations,
such as the determination of the $10^2$ central eigenpairs of a $10^9$-dimensional matrix presented here.
To our knowledge presently no other method is capable of computations at this scale.
These promising results indicate that ChebFD can become a viable tool for electronic structure calculations in quantum chemistry and physics, and for research into modern topological materials.

The present ChebFD scheme can be refined along several directions.
First, the construction of better polynomial filters could improve convergence. 
Second, adaptive schemes that adjust the polynomial degree and the number of search vectors according to the DOS, which can be computed with increasing accuracy in the course of the iteration, might prove useful.
Third, Ritz values should be replaced by harmonic Ritz values in the computation, together with an improved convergence criterion and locking of converged eigenpairs.
Note, however, that our theoretical analysis of ChebFD shows that none of these improvements is truly crucial for the (fast) convergence of polynomial filtering: The most important aspect is the ``overpopulation'' of the target space with search vectors, i.e., the use of a sufficiently large search space.

A more radical alternative to ChebFD is the replacement of polynomial filters by rational functions as, e.g., in the FEAST and CIRR methods.
From the theoretical point of view, rational functions always lead to better convergence.
However, evaluation of rational functions of sparse matrices is non-trivial. Standard iterative solvers require too many additional spMVMs to be competitive with polynomial filters, unless they can exploit additional structural or spectral properties of the matrix~\cite{2014-GalgonKraemerEtAl-OnTheParallelIterative-PREPRINT}. At present it is not clear whether non-polynomial filter functions can succeed in large-scale problems of the size reported here. This question is currently under investigation.

The ChebFD algorithm implemented and investigated in the present work should be of particular interest to physicists and chemists who increasingly require high-performance tools for the computation of interior eigenvalues in present and future research applications. 
The \ghost\ library, the ChebFD code, and generator routines for all test cases are publicly available from the ESSEX repository\footnote{\url{https://bitbucket.org/essex}}.

\section{Acknowledgments}
The research reported here was funded by Deutsche Forschungsgemeinschaft through priority programmes 1459 ``Graphene'' and 1648 ``Software for Exa\-scale Computing''.
The authors gratefully acknowledge the Gauss Centre for Supercomputing e.V. (\url{http://www.gauss-centre.eu}) for funding this work by providing computing time on the GCS Supercomputer SuperMUC at Leibniz Supercomputing Centre (LRZ, \url{http://www.lrz.de}) through project pr84pi.
The \ghost-based implementation of ChebFD will become available as part of the sparse solver library prepared in the ESSEX project\footnote{\url{http://blogs.fau.de/essex}}.

\end{document}